\documentclass{amsart}
\usepackage{epsf}
\usepackage{epsfig}
\newcommand{\Img}[2]{\epsfxsize#1truecm\epsfbox{#2}}
\def\C{{\mathbb{C}}}
\def\Z{\mathbb{Z}}
\def\R{\mathbb{R}}
\def\M{{\mathcal M}}
\def\sm{\operatorname{sm}}
\def\cm{\operatorname{cm}}
\def\smh{\operatorname{smh}}
\def\cmh{\operatorname{cmh}}
\def\ds{\displaystyle}
\def\Inv{\operatorname{Inv}}
\def\F{{\bf F}}

\def\tode{\,\mathop{\longrightarrow}^\delta\,}
\def\Pr{\mathbb{P}}
\def\cal{\mathcal}
\def\Tree{\operatorname{Tree}}
\def\b{{\sf x}}
\def\w{{\sf y}}
\def\x{{\sf x}}
\def\y{{\sf y}}
\def\lora{\longrightarrow}

\def\Lap{{\cal L}}
\def\frak{\mathfrak}
\newcommand{\der}[1]{\frac{\partial}{\partial #1}}
\renewcommand{\arraystretch}{1.25}

\newtheorem{theorem}{Theorem}
\newtheorem{proposition}{Proposition}
\newtheorem{definition}{Definition}

\newcounter{noteno}

\newenvironment{Note}%
	{\refstepcounter{noteno}%
	\begin{small}
	\medbreak\par\noindent{{\bf Note~\thenoteno}.}}%
	{\hfill{$\Box$}\end{small}\par\medbreak}

\begin{document}

\title[A Combinatorial Excursion]{\bf The Fermat cubic, elliptic functions,
continued fractions, and a combinatorial excursion}
\author{Eric van Fossen Conrad and Philippe Flajolet}
\address{E.C.: Department of Mathematics, The Ohio State University,
231 W 18th Avenue, 
Columbus, OH 43210 USA }
\email{econrad@math.ohio-state.edu}
\address{P.F.: Algorithms Project, INRIA Rocquencourt, F-78153
(France)}
\email{Philippe.Flajolet@inria.fr}
\dedicatory{Kindly dedicated to G\'erard\,$\cdots$Xavier Viennot on
the occasion of his sixtieth birthday.} 
\date{July 9, 2005}

\thanks{{\bf Submitted to \emph{S\'eminaire Lotharingien de
Combinatoire}  (Journal).} \\
Presented at the 54th \emph{S\'eminaire Lotharingien de
Combinatoire} (Conference), Lucelle, April 2005. The work of E. Conrad 
was supported in part by National Security Agency grants
MDA 904-93-H-3032, MDA 904-97-1-0019 and MDA 904-99-1-0003, and by
National Science Foundation grant DMS-0100288. The work of P. Flajolet was
supported in part by a grant from ACI-NIM (New Interfaces of Mathematics).}

\begin{abstract}
Elliptic functions considered by Dixon in the nineteenth century and
related to Fermat's cubic, $x^3+y^3=1$, lead to a new set
of continued fraction expansions with sextic numerators
and cubic denominators. The functions and the fractions are 
pregnant with interesting combinatorics, including a special P\'olya urn, a
continuous-time branching process of the Yule type, as well as permutations 
satisfying various constraints that involve either parity of levels of
elements or a repetitive pattern of order three. The combinatorial
models are 
related to but different from models of elliptic functions earlier
introduced by Viennot, Flajolet, Dumont, and Fran{\c c}on. 
\end{abstract}

\maketitle

In 1978, Ap\'ery announced an amazing
discovery: \emph{``$\zeta(3)\equiv \sum 1/n^3$ is irrational''}. 
This represents a great piece of Eulerian mathematics
of which van der Poorten has written a particularly
vivid account in~\cite{Poorten79}.
At the time of Ap\'ery's result, nothing was known about
the arithmetic nature of the zeta 
values at odd integers, and not unnaturally his theorem 
triggered interest in a whole range of problems that are now recognized
to relate to
much ``deep'' mathematics~\cite{KoZa01,Rivoal02}.
Ap\'ery's original irrationality proof crucially
depends on a continued fraction representation
of~$\zeta(3)$. To wit:
\begin{equation}\label{zeta3}\renewcommand{\arraystretch}{1.4}
\begin{array}{l}
\zeta(3)= \ds\cfrac{6}{\varpi(0)-\cfrac{1^6}
	{\varpi(1)-\cfrac{2^6}{\varpi(2)-\cfrac{3^6}{{\ddots}}}}},
\\
\hbox{where}\qquad \varpi(n)=(2n+1)(17n(n+1)+5).
\end{array}
\end{equation}
(This is not of a form usually considered 
by number theorists.) What is of special interest
to us is that the $n$th stage of the fraction involves the
\emph{sextic} numerator $n^6$, while the corresponding numerator is a 
\emph{cubic}
polynomial in~$n$. Mention must also be made at this stage of a
fraction due to Stieltjes (to be later rediscovered
and  extended by Ramanujan~\cite[Ch.~12]{Berndt89}), namely,
\begin{equation}\label{psi2}\renewcommand{\arraystretch}{1.4}
\begin{array}{l}
\ds \sum_{k=1}^\infty \frac{1}{(k+z)^3}~=~ \cfrac{1}{\sigma(0)-\cfrac{1^6}
	{\sigma(1)-\cfrac{2^6}{\sigma(2)-
\cfrac{3^6}{\ddots}}}},
\\
\hbox{where}\qquad \sigma(n)=(2n+1)(n(n+1)+2z(z+1)+1).
\end{array}
\end{equation}
Unfortunately that one  seems to have no useful arithmetic content.

Explicit continued fraction expansions of special functions are
really \emph{rare}~\cite{Perron54,Wall48}. Amongst the very few known,
most involve  
numerators and denominators which, at depth~$n$ in the continued fraction,
depend rationally on~$n$ in a manner that is at most quartic.
In this context, the Stieltjes-Ramanujan-Ap\'ery fractions stand out.
It then came as a surprise that certain
functions related to Dixon's nineteenth century 
parametrization of Fermat's cubic $X^3+Y^3=1$ lead
to continued fractions that precisely share with~(\ref{zeta3})
and~(\ref{psi2}) the cubic--sextic dependency of the their
coefficients on the depth~$n$, for instance,
\begin{equation}\label{conrad-cf}
\begin{array}{l}
 \ds \int_0^\infty \sm(u) e^{-u/x}\, du =
\cfrac{x^2}{1+b_0x^3-\cfrac{1\cdot2^2\cdot3^2\cdot 4\, x^6}
	{1+b_1x^3-\cfrac{4\cdot5^2\cdot6^2\cdot7\, x^6}
	{1+b_2x^3-\cfrac{7\cdot8^2\cdot9^2\cdot10\,
x^6}{{\ddots}}}}}\, ,
\\
\ds
\hbox{where}\qquad 
b_n=2(3n+1)((3n+1)^2+1).
\end{array}
\end{equation}
There, the function $\sm$ is in essence the inverse of a
${}_2F_1$--hypergeometric of type $F[\frac13,\frac23;\frac43]$;
see below for proper definitions.
This discovery, accompanied by several related
continued fractions, 
was first reported in Conrad's PhD thesis~\cite{Conrad02}
defended in~2002. It startled the second author
with a long standing interest in continued
fractions~\cite{Flajolet80b,Flajolet82,FlFr89}, 
when he discovered from reading Conrad's thesis
in early 2005, that certain elliptic functions 
could precisely lead to a cubic--sextic fraction. 
This paper describes our ensuing exchanges. 
We propose to show
that there is an interesting orbit of ideas and results surrounding the Fermat
cubic, Dixon's elliptic functions,  Conrad's fractions, 
as well as certain urn models
of probability theory, Flajolet's
combinatorial theory of continued fractions, and finally,
the elementary combinatorics of permutations.

\medskip
\noindent
{\bf Plan of the paper.}
The Dixonian elliptic functions, ``$\sm$'' and ``$\cm$,
are introduced in Section~\ref{fermat-sec}.
Their basic properties are derived
from a fundamental nonlinear differential system that they satisfy.
It this way, one can prove simply that they parametrize the Fermat
cubic and at the same time admit of representations as inverses of
hypergeometric functions. Section~\ref{conrad-sec} presents the
complete proof 
of six continued fractions of the Jacobi type and three fractions of 
the Stieltjes type that are associated (via a Laplace transform)
to Dixonian functions---this is
the first appearance in print of results from Conrad's PhD
thesis~\cite{Conrad02}. Next, in Section~\ref{balls-sec}, we prove that the 
combinatorics of the nonlinear differential system defining $\sm,\cm$ 
is isomorphic to the stochastic evolution
of a special process,
which is an urn of the P\'olya type. 
As a consequence, this urn, together with its continuous-time analogue (a
two-particle version of a classical binary branching process) can be solved
analytically in terms of Dixonian functions. Conversely, this isomorphism 
constitutes a
first interpretation of the coefficients of $\sm,\cm$ phrased in terms
of combinatorial objects that are \emph{urn histories}. Our second
combinatorial interpretation is in terms of \emph{permutations}
and it involves \emph{peaks, valleys, double rises, and double falls} as well as
the parity of levels of nodes in an associated tree representation;
see Section~\ref{firstperm-sec}.
(Several combinatorial models of elliptic 
functions due to Dumont, Flajolet, Fran{\c c}on,
and Viennot~\cite{Dumont79,Dumont81,Flajolet80b,FlFr89,Viennot80}
are otherwise known to involve parity-constrained permutations.)
Next, the continued fraction expansions relative to Dixonian
functions can be read combinatorially through the glasses of 
a theorem of Flajolet~\cite{Flajolet80b}
and a bijection due to Fran{\c c}on and Viennot~\cite{FrVi79},
relative to systems of weighted lattice paths and continued fractions.
This is done in Section~\ref{secondperm-sec} which presents our third
combinatorial model of Dixonian functions: the coefficients of
$\sm,\cm$ are shown to enumerate certain types of \emph{permutations involving a
repetitive pattern} of order~three. 
Finally, Section~\ref{connec-sec}  briefly summarizes a few other
works where Dixonian functions make an appearance.

\smallskip\noindent
{\bf Warning.} This paper corresponds to an invited lecture at the
\emph{Viennotfest} and, as
such, its style is often informal. Given the scarcity of the
literature relative to Dixonian functions, we have attempted to provide pointers to all of the relevant
works available to us. Thus, our article attempts to kill three birds
with one stone, namely be a tribute to Viennot, survey the area,
and present original results.

\smallskip\noindent
{\bf Dedication.}
This paper is kindly dedicated to G\'erard$\ldots$Xavier Viennot on the
occasion of his sixtieth birthday, celebrated at Lucelle in April 2005. His works in lattice path
enumeration, bijective combinatorics, and the combinatorics of
Jacobian elliptic functions have inspired us throughout
the present work.

\section{The Fermat curves, the circle, and the cubic}\label{fermat-sec}
The Fermat curve $\F_m$ is the complex algebraic curve defined by the
equation
\[
X^m+Y^m=1.
\]
(Fermat-Wiles asserts that this curve has no 
nontrivial rational point
as soon as $m\ge3$.)

Let's start with the innocuous looking $\F_2$,
that is, the \emph{circle}. Of interest for this discussion is
the fact that the circle can be parametrized by trigonometric
functions. Consider the two functions from $\C$ to $\C$ defined by
the \emph{linear} differential system,
\begin{equation}\label{circ1}
s'=c,\quad c'=-s,\qquad
\hbox{with initial conditions}\quad
s(0)=0,\quad c(0)=1.
\end{equation}
Then the transcendental functions $s,c$ do parametrize the circle, since
\[
s(z)^2+c(z)^2=1,
\]
as is verified immediately from the differential system, which implies
$(s^2+c^2)'=0$.
At this point one can switch to conventional notations and set
\[
s(z)\equiv \sin z,\qquad
c(z)\equiv \cos(z).\]
It is of interest to note that these functions are also obtained by
inversion from an \emph{Abelian integral}\footnote{
Given an algebraic curve $P(z,y)=0$, an Abelian integral is any integral 
$\int R(z,y)\, dz$, where $R$ is a rational function.}
 on the $\F_2$ curve:
\[
\int_0^{\sin z}\frac{dt}{\sqrt{1-t^2}},\qquad
\cos(z)=\sqrt{1-\sin(z)^2}
\]
For combinatorialists, it is of special interest to note that the related 
functions
\[
\tan z=\frac{\sin z}{\cos z}, \qquad
\sec z= \frac{1}{\cos z},\]
enumerate a special class of permutations, the alternating
(also known as ``up-and-down'' or ``zig-zag'') ones. This last fact
is a classic result of combinatorial analysis discovered by D\'esir\'e Andr\'e around 1880.

\subsection{The Fermat cubic and its Dixonian parametrization}
Next to the circle, in order of complexity, comes the Fermat cubic $\F_3$. 
Things should be less elementary since the Fermat curve has
(topological) genus~1, but this very fact points to strong
connections with  elliptic functions. 

\begin{figure}
\begin{center}
\hbox{\Img{5}{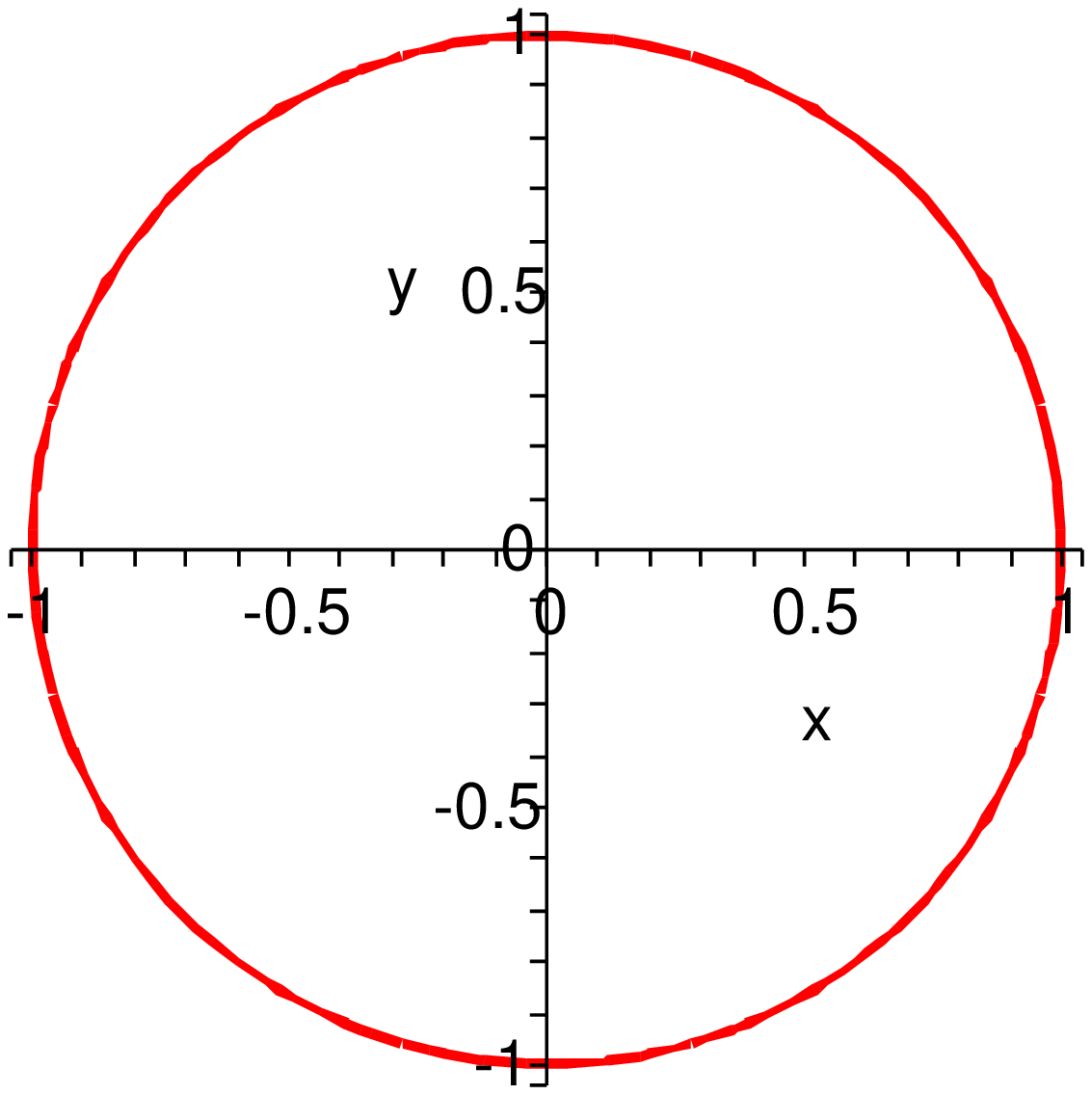}\qquad \Img{5}{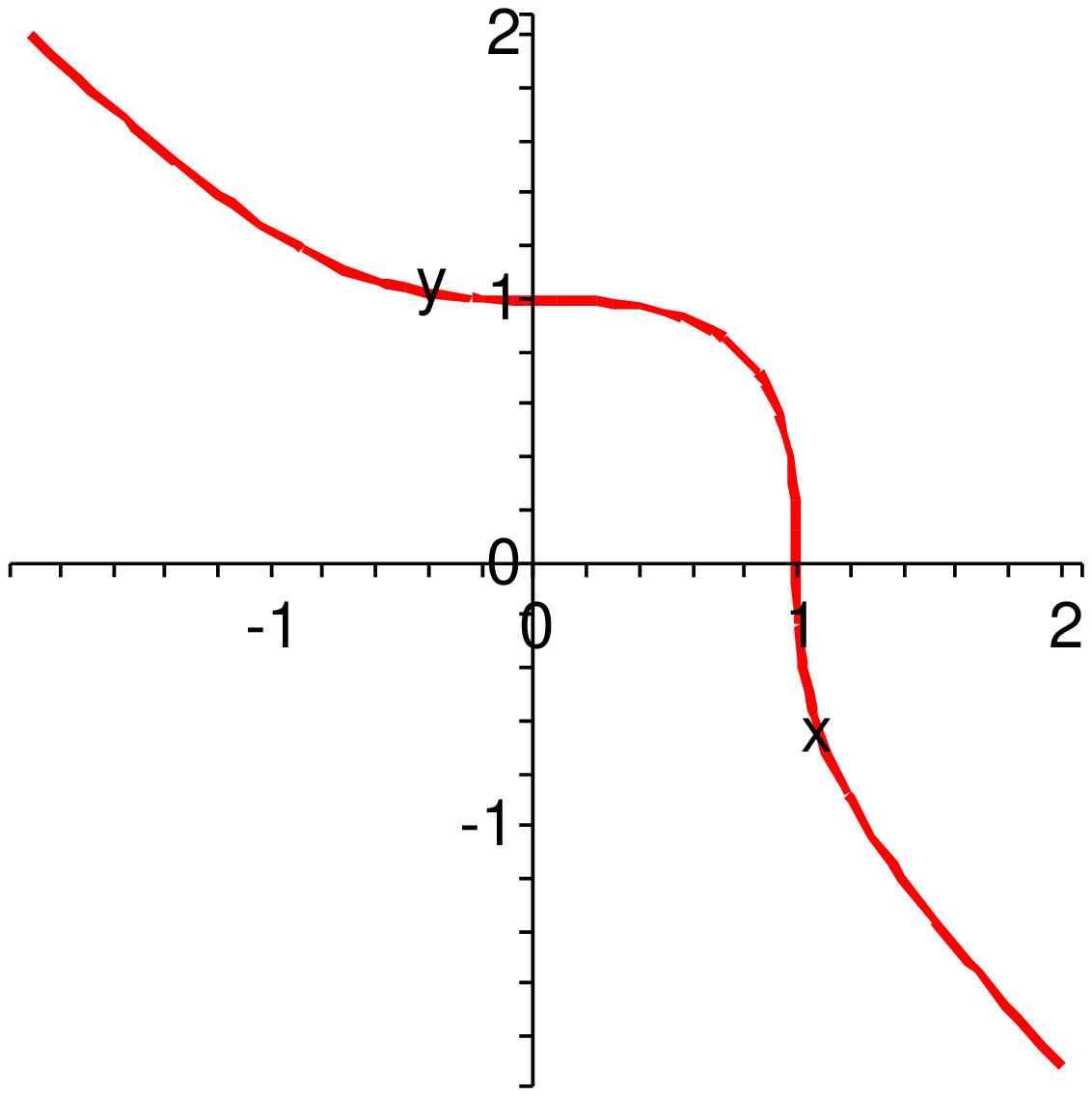}}
\end{center}
\caption{\label{fermat-fig}
The two Fermat curves $\F_2$ and $\F_3$.}
\end{figure}

The starting point is a clever generalization of~(\ref{circ1}).
Consider now the \emph{nonlinear} differential system
\begin{equation}\label{conrad1}
s'=c^2,\qquad c'=-s^2,
\end{equation}
with initial conditions $s(0)=0$, $c(0)=1$. 
These functions are analytic about the origin, a fact resulting from
the standard existence theorem for solutions of ordinary differential
equations. 
Then, a one line
calculation  similar to the trig function case shows that
\[
s(z)^3+c(z)^3=1,\]
since
\[
\left(s^3+c^3\right)'=3s^2c^2-3c^2s^2=0.\]
Consequently, the pair $\langle s(z),c(z)\rangle$ parametrizes the Fermat
curve $\F_3$, at least \emph{locally} near the point~$(0,1)$.
The basic properties of these functions have been elicited by 
 Alfred Cardew Dixon (1865--1936) in a long paper~\cite{Dixon90}.
Dixon established that the functions are meromorphic in the whole
of the complex plane and doubly periodic (that is, elliptic), hence
they provide a \emph{global} parametrization of the Fermat cubic.

\begin{figure}
\Img{5.5}{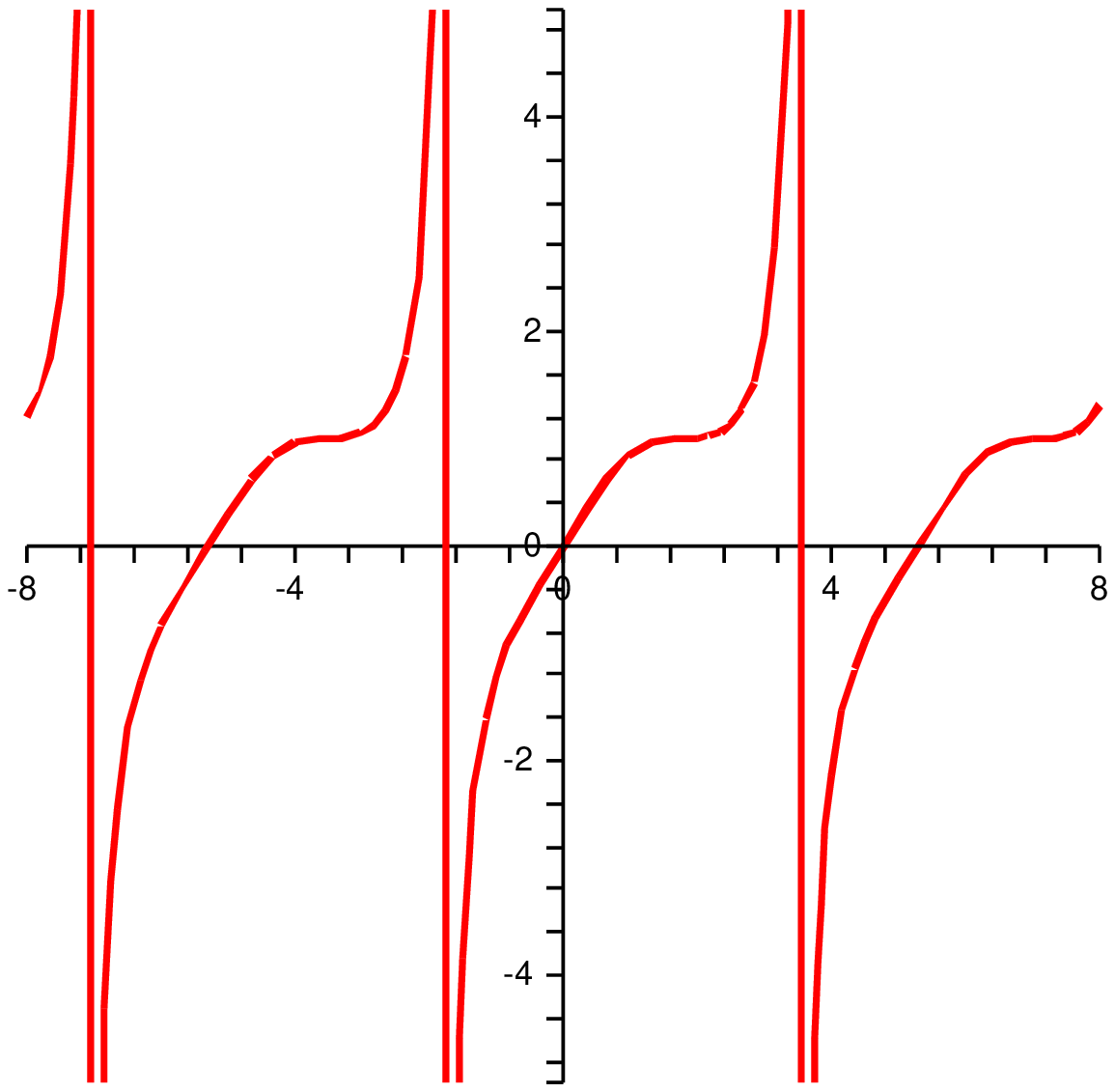}\hspace*{0.5truecm}\Img{5.5}{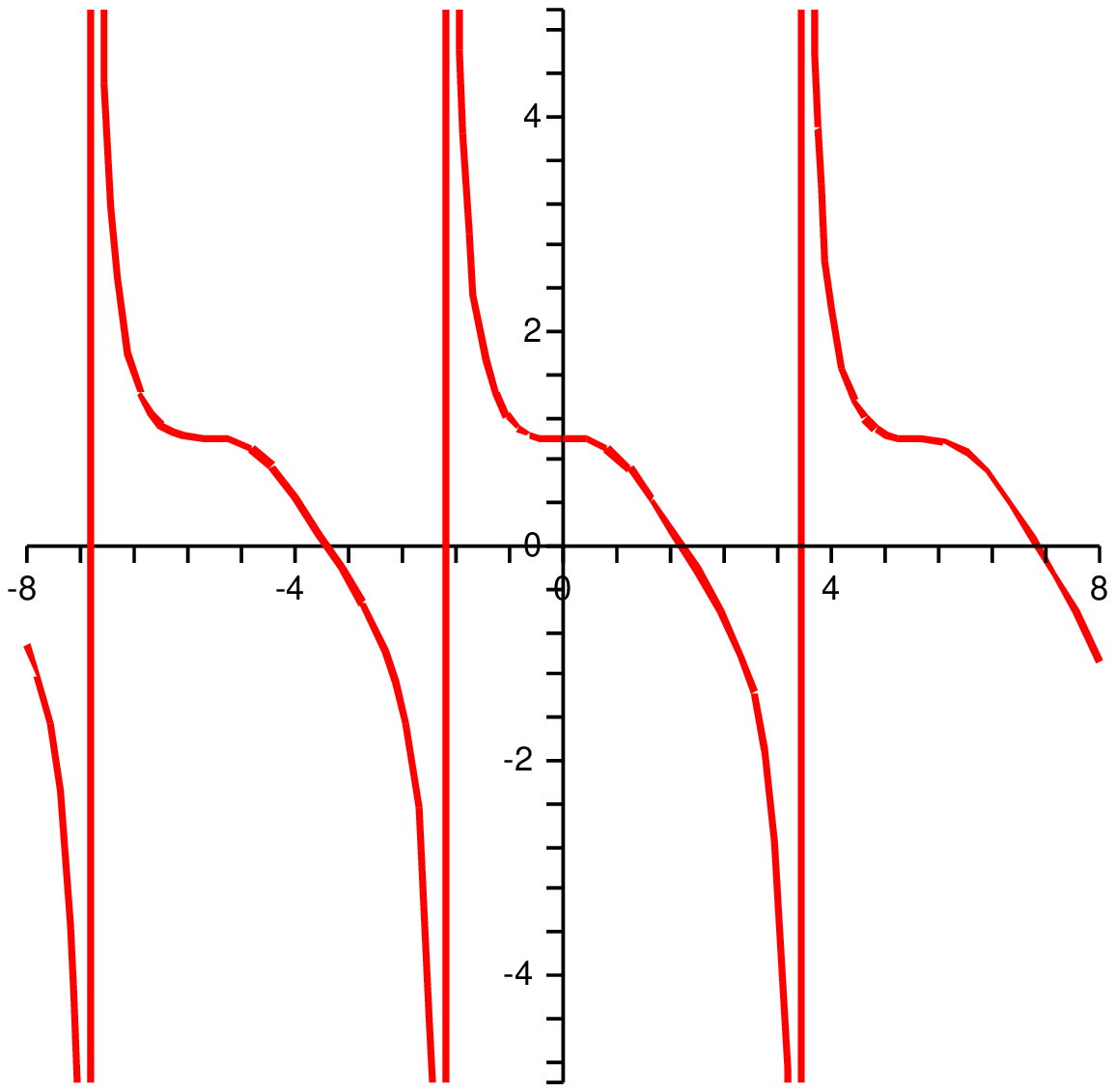}
\caption{\label{thenumplots-fig}
Plots of $\sm(z)$ [left] and $\cm(z)$
[right] for $z\in\R$.}
\end{figure}

From now on, we'll give the $s,c$ functions
the name introduced by Dixon (see Figure~\ref{thenumplots-fig} for a
rendering)
and set
\[
\sm(z)\equiv s(z), \qquad \cm(z)\equiv c(z).\]
Their Taylor expansions  at~0  (not currently found in Sloane's 
\emph{Encyclopedia of Integer Sequences (EIS)}~\cite{Sloane00})
start as follows:
\begin{equation}\label{sctayl}
\left\{\begin{array}{lll}
\sm(z)&=& \ds z-4\frac{z^4}{4!}+160
\frac{z^7}{7!}-20800\frac{z^{10}}{10!}+6476800\frac{z^{13}}{13!}-\cdots
\\
\cm(z)&=& \ds
1-2\frac{z^3}{3!}+40\frac{z^6}{6!}-3680\frac{z^9}{9!}+
8880000 \frac{z^{12}}{12!}-\cdots\,.
\end{array}\right.
\end{equation}
In summary, we shall operate with the new notations $\sm,\cm$
and with the defining system:
\begin{equation}\label{dixon0}
\sm'(z)=\cm^2(z),\quad \cm'(z)=-\sm^2(z);\qquad 
\sm(0)=0,\quad \cm(0)=1.
\end{equation}

\subsection{A hypergeometric connection.}\label{hyper-sec}
At this point, it is worth noting that one can easily make the
$s\equiv\sm$ and $c\equiv\cm$ functions somehow ``explicit'' 
Start from the defining system 
(System~(\ref{dixon0}) abbreviated here as $(\Sigma)$) and 
apply differentiation ($\partial$):
\begin{equation}\label{sol1}
s'=c^2 \quad\mathop{\Longrightarrow}^{\partial}\quad
s''=2cc' \quad\mathop{\Longrightarrow}^{\Sigma}\quad
s''=-2cs^2 \quad\mathop{\Longrightarrow}^{\Sigma}\quad
s''=-2s^2 \sqrt{s'}.
\end{equation}
Then ``cleverly'' multiply by $\sqrt{s'}$ to get, via integration ($\int$),
\begin{equation}\label{sol2}
s''\sqrt{s'}=-2s^2s' \quad\mathop{\Longrightarrow}^{\int}\quad
\frac23 (s')^{3/2}=-\frac23 s^3+K,
\end{equation}
for some integration constant $K$ which must be equal to~$\frac23$, given the
initial conditions. This proves in
two lines(!) that $\sm$ is the inverse of an
integral,
\begin{equation}\label{abel3}
\int_0^{\sm(z)} \frac{dt}{(1-t^3)^{2/3}}=z,
\end{equation}
this integral being at the same time an incomplete Beta integral and
an Abelian integral over the Fermat curve ($\int \frac{dy}{y^2}$).
By the same devices, it is seen that the function~$\cm(z)$ also satisfies
\begin{equation}\label{cmint}
z=\int_{\cm z}^1 \frac{dt}{(1-t^3)^{2/3}},
\end{equation}
and is thus itself the inverse of an Abelian integral.

Via expansion and term-wise integration,
the latest finding~(\ref{abel3}) can then be rewritten in terms of the
classical hypergeometric function,
\[
{}_2F_1[\alpha,\beta,\gamma;z]:=1+\frac{\alpha\cdot\beta}{\gamma}\frac{z}{1!}
+\frac{\alpha(\alpha+1)\cdot\beta(\beta+1)}{\gamma(\gamma+1)}\frac{z^2}{2!}
+\cdots\,.
\]
with special rational values of the parameters. Letting $\Inv(f)$ denote the
inverse of~$f$ with respect to composition (i.e., $\Inv(f)=g$ if
$f\circ g=g \circ f =Id$) we state:

\begin{proposition}[Dixon~\cite{Dixon90}]\label{prop-hyper}
The function $\sm$ is the inverse of an Abelian integral over~$\F_3$
and equivalently the inverse of a ${}_2F_{1}$:
\[
\sm (z)= \Inv \int_0^z \frac{dt}{(1-t^3)^{2/3}}
=\Inv z\cdot {}_2F_1\left[\frac{1}{3},\frac{2}{3},\frac{4}{3} ;z^3\right].
\]
The function $\cm$ is then defined near~$0$ by $\cm(z)=\root 3 \of {1-\sm^3
(z)}$, or alternatively by~\eqref{cmint}.
\end{proposition}

The analogy with the sine function is striking.
Of course, Proposition~\ref{prop-hyper} is not new and all this is
related to extremely classical material.
Dixon~\cite{Dixon90} discusses the implicit integral representations~(\ref{abel3}),
~(\ref{cmint}) 
and writes concerning the prehistory of his $\sm,\cm$:
\begin{quote}\small
	The only direct references that I have come across elsewhere
	are certain passages in the lectures of Professor Cayley and
	Mr. Forsyth where the integral $\int\frac{dx}{(1-x^3)^{2/3}}$
	was used to illustrate Abel's Theorem, in the treatises of Legendre,
	and Briot and Bouquet, and again in Professor Cayley's
	lectures and elsewhere where it is shewn how to turn the integral
	into elliptic form, and lastly at the end of Bobek's
	{\em Einleitung in die Theorie der elliptischen Functionen} where
	expressions are found for the coordinates of any point on the
	above cubic.
\end{quote}
It is fascinating to be able to develop  a fair amount of
the theory from the differential equation~(\ref{dixon0}),
using only first principles of analysis. (See Dixon's article as well as
Section~\ref{connec-sec} of the present paper for more.)

\begin{Note} \emph{ Lundberg's hypergoniometric functions.}
The question of higher degree generalizations of the differential 
system satisfied
by $\sm,\cm$ is a natural one. Indeed, the system
$s'=c^{p-1},c'=-s^{p-1}$ parametrizes locally near $(0,1)$ the Fermat
curve $X^p+Y^p=1$. (This parametrization ceases to be a global one,
however, since $\F_p$ has genus $g=(p-1)(p-2)/2$, so that $g\ge3$ as
soon as $p\ge4$; see for instance~\cite{Lang82}.
Thus, only the case $p=3$ leads directly to elliptic
functions.) The corresponding functions are
still locally inverses of Abelian integrals over the Fermat curve, which is
verified by calculations similar to Equations~(\ref{sol1})
and~(\ref{sol2}). 
In fact, a rather unknown high
school teacher
from Sweden, Erik Lundberg, 
developed in 1879  an elementary
theory (see~\cite{Lundberg79}) of what he called ``hypergoniometric functions'',
including a \emph{sinualis} and a \emph{cosinualis} that are indexed
by rational numbers. Interest in these questions was recently
rekindled by an insightful article of Lindqvist and Peetre
published in the \emph{American Mathematical
Monthly}~\cite{LiPe01}. The authors
discuss various connections to elliptic functions, in particular,
the reduction of Fermat's cubic to its Weierstra{\ss} normal form by
elementary manipulations.
See also a problem~\cite{LiPe01b}  in the same issue of the 
\emph{Monthly} posed by these
authors together with a remark due to Jon Borwein to the effect that
\[
\sm\left(\frac{\pi_3}{3}-x\right)=\cm(x).
\]
There $\pi_3$ is a fundamental constant of Dixonian functions:
\begin{equation}\label{pi3def}
\pi_3=3\int_0^1
\frac{dt}{(1-t^3)^{2/3}}=B\left(\frac13,\frac13\right)\equiv
\frac{\Gamma\left(\frac13\right)^2}{\Gamma\left(\frac23\right)}
=\frac{\sqrt{3}}{2\pi}\Gamma\left(\frac13\right)^3\doteq
5.29991\,62508.
\end{equation}
Many more properties of this type are to be found in Dixon's paper.
In particular, one has, 
\[
\sm(\pi_3+u)=\sm(u),\qquad \cm(\pi_3+u)=\cm(u),
\]
that is, $\pi_3$ is a real period. The complete lattice of periods of
$\sm,\cm$ is
\[
\Z\pi_3\oplus \Z\pi_3\omega,\qquad
\omega:=e^{2i\pi/3},\]
which is consistent with the rotational
invariance: $\sm(\omega
u)=\omega\sm(u)$ and $\cm(\omega u)=\cm(u)$. It thus corresponds to 
the hexagonal lattice displayed on Figure~\ref{sixlat-fig} 
and to case $C$ of the urn models evoked in
Subsection~\ref{sixlat-sec}.
%
\end{Note}

\section{Some startling fractions}\label{conrad-sec}

In 1907, L.\ J.\ Rogers~\cite{Rogers07} devised two methods to
obtain continued fraction expansions of
Laplace transforms of the Jacobian elliptic functions $sn,cn$
(see for instance~\cite{WhWa27} for the basic theory of these functions).
In his first method, he resorted to integration by parts.  His
second method involves 
a general ``addition theorem''~\cite{Perron54,Wall48};
it was
to some extent a rediscovery 
of a technique introduced by T.\ J.\ Stieltjes in
\cite{Stieltjes89} that relies on diagonalization of certain infinite
quadratic forms.              
Stieltjes and Rogers found altogether three families of continued fractions
relative to $sn,cn$.  (A fourth family, implicit in the work of Rogers,
was discovered almost
a century later by Ismail and Masson~\cite{IsMa99}.) 
Such expansions turn out to be useful: S.\ Milne in
\cite{Milne02} obtained additional expansions
implying explicit Hankel determinant evaluations,
which enabled him to
prove some deep results about sums of squares and sums of triangular
numbers.   Milne's results in particular include  exact
explicit infinite families of identities expressing
the number of ways to write an integer as the sum of $4N^2$ or $4N(N+1)$
squares of integers, where $N$ is an arbitrary positive integer.

In his PhD thesis defended in 2002, Conrad~\cite[Ch.~3]{Conrad02}
applied the integration by parts method
to develop completely
new continued fractions arising from Dixon's elliptic functions. 
These fractions fall into six families
which are naturally grouped as two sets of three.
The underlying symmetries suggest that these classifications are
fundamental and complete.

In what follows, we make use of the \emph{Laplace transform} classically 
defined by
\[
\Lap(f,s)=\int_0^\infty f(u)e^{-su}\, du.
\]
For our purposes, it turns out to be convenient to set $s=x^{-1}$.
In that case, one has
\[
x^{-1}\Lap\left(u^\nu,x^{-1}\right)= \nu! x^n,\]
which means that the Laplace transform 
formally maps exponential generating functions (EGFs) to ordinary
generating functions (OGFs):
\begin{equation}\label{formlap}
x^{-1}\Lap\left(\sum_{\nu\ge0} a_\nu \frac{u^\nu}{\nu!},x^{-1}\right)=
\sum_{\nu\ge0} a_\nu {x^\nu}.
\end{equation}
There are a number of conditions ensuring the analytic or
asymptotic validity of this last equation. In what follows, we shall
make use of the integral notation for the Laplace transformation,
but only use it as a convenient way to represent the \emph{formal transformation} from 
EGFs to OGFs in~(\ref{formlap}).

The continued fractions that we derive are of two types~\cite{Perron54,Wall48}. The first
type, called a $J$-fraction ($J$ stands for Jacobi), associates to a
formal power series $f(z)$ a fraction whose denominators are linear
functions of the variable~$z$ and numerators are quadratic monomials:
\[
f(z)=\cfrac{1}{1-c_0z-\cfrac{b_1z^2}{1-c_1z-\cfrac{b_2z^2}{\ddots}}}\,.
\]
(Such fractions are the ones naturally associated to orthogonal
polynomials~\cite{Chihara78}.) 
The second type\footnote{%
	An alternative name for $S$-fraction is ``regular
$C$-fraction''. What we called $J$-fraction is also known as
an ``associated continued fraction''.
}, called in this paper an $S$-fraction ($S$ stands for Stieltjes) 
has denominators that are the constant~1 and numerators that are
monomials of the first degree:
\[
f(z)=\cfrac{1}{1-\cfrac{d_1z}{1-\cfrac{d_2z}{\ddots}}}\,.
\]
An $S$-fraction can always be contracted into a~$J$-fraction, 
with the corresponding formul{\ae} being explicit, but
not conversely. Accordingly, from
the point of view of the theory of special functions, an explicit
$S$-fraction expansion should be regarded as a stronger form than
its $J$-fraction counterparts. 

\subsection{$J$-fractions for the Dixon functions.}
We introduce three classes of formal Laplace transforms:
\begin{equation}\label{SCD}
\renewcommand{\arraystretch}{1.9}
\begin{array}{lll}
S_n &:=& \ds \int_0^\infty \sm(u)^n \, e^{-u/x} du \\
	C_n  &:=& \ds \int_0^\infty \sm(u)^n \cm(u) \, e^{-u/x} du \\
	D_n &:=& \ds \int_0^\infty \sm(u)^n \cm(u)^2 \, e^{-u/x} du
\end{array}
\end{equation}
(Since $\sm^3+\cm^3=1$, we can reduce any polynomial in $\sm$ and $\cm$
to one which is 
at most of degree 2 in $\cm$.  We let powers of $\sm$ grow since $\sm u$
vanishes at $u=0$.)
After integrating by parts and reducing to canonical form using the Fermat
cubic parametrization, we obtain the following recurrences:
\begin{equation}
\label{Dixon:rec:first}
\begin{array}{lll}
	S_0 = x, 
	& S_n  = n x D_{n-1} &  (n>0)  \\
	C_0 = x - x S_2, 
	& C_n = n x S_{n-1}  - (n+1) x S_{n+2} & (n>0)  \\
	D_0  = x - 2 x C_2, 
	& D_n  = n x C_{n-1}  - (n+2) x C_{n+2} & (n>0)
\end{array}
\end{equation}
These recurrences will provide six $J$-fractions and three $S$-fractions.

For Laplace transforms of powers of $\sm(u)$, the
recurrences~(\ref{Dixon:rec:first}) yield the fundamental
relation,
\begin{equation}\label{smn}
	\frac{S_n}{S_{n-3}} =
		\frac{n(n-1)(n-2)x^3}{
			1 + 2n(n^2+1)x^3 - n(n+1)(n+2)x^3
\frac{S_{n+3}}{S_n} },
\end{equation}
subject to the initial conditions
\[
	S_1 = \frac{x^2}{1 + 4x^3 - 6x^3 \frac{S_4}{S_1}}, \quad
	S_2 = \frac{2x^3}{1 + 20x^3 - 24x^3 \frac{S_5}{S_2}}, \quad
	S_3 = \frac{6x^4}{1 + 60x^3 - 60x^3 \frac{S_6}{S_3}}.
\]
Repeated use of the relation~(\ref{smn}) starting from the initial
conditions then leads to three continued fraction expansions
relative to each of $S_1,S_2,S_3$.

An entirely similar process applies to $C_n$. We find
\[
 	\frac{C_n}{C_{n-3}} =
		\frac{n(n-1)(n-2)x^3}{
			1 + \left((n-1)n^2 +(n+1)^2(n+2)\right) x^3
				- (n+1)(n+2)(n+3)x^3
\frac{C_{n+3}}{C_n} }
,\]
subject to
\[
	C_0 = \frac{x}{1 + 2x^3 - 6x^3 \frac{C_3}{C_0}}, \quad
	C_1 = \frac{x^2}{1 + 12x^3 - 24x^3 \frac{C_4}{C_1}}, \quad
	C_2 = \frac{2x^3}{1 + 40x^3 - 60x^3 \frac{C_5}{C_2}}
.\]

For $D_n$,
we could handle things
in a similar manner, but it turns out that there is a simpler
procedure available 
to us.  Differentiating powers of $\sm(u)$ with respect to $u$, we have
\[
	\frac{d}{du} \sm(u)^{n+1} = (n+1) \, \sm(u)^n \, \cm^2(u)
\text{.}
\]
Given the effect of Laplace transforms on derivatives, 
it follows that 
\[
	\mathcal{L} ( \sm(u)^n \, \cm(u)^2, x^{-1} )
		= \frac{1}{(n+1)x} \mathcal{L} ( \sm^{n+1}(u), x^{-1} )
\text{.}
\]
We thus obtain no substantially new $J$-fractions from $D_n$
by this process.

\begin{theorem}[Conrad~\cite{Conrad02}] The formal Laplace transform
of the functions $\sm^n$ $(n=1,2,3)$ and  $\sm^n\cdot\cm$ $(n=0,1,2)$
 have explicit 
$J$-fraction expansions with cubic denominators and sextic numerators
as given in Figure~\ref{sixJ-fig}.
\end{theorem}

\begin{figure}[p]\small
\fbox{\begin{minipage}{\textwidth}
\begin{align*}
	\int_0^\infty \sm u \, e^{-u/x} \, du &=
		\cfrac{x^2}{1+4x^3-
			\cfrac{a_1 x^6}{ 1 + b_1 x^3 -
			\cfrac{a_2 x^6}{ 1 + b_2 x^3 - \ddots}}}
\intertext{where\quad
	$a_n = (3n-2)(3n-1)^2 (3n)^2 (3n+1)$,\quad
	$b_n = 2(3n+1) \left( (3n+1)^2 + 1 \right)$.}
	\int_0^\infty \sm^2(u) \, e^{-u/x} \, du &=
		\cfrac{2 x^3}{1+20x^3-
			\cfrac{a_1 x^6}{ 1 + b_1 x^3 -
			\cfrac{a_2 x^6}{ 1 + b_2 x^3 - \ddots}}}
\intertext{where\quad
	$a_n = (3n-1)(3n)^2 (3n+1)^2 (3n+2)$,\quad
	$b_n = 2(3n+2) \left( (3n+2)^2 + 1 \right)$.}
	\int_0^\infty \sm^3(u) \, e^{-u/x} \, du &=
		\frac{6x^4}{1+60x^3+
			\cfrac{a_1 x^6}{ 1 + b_1 x^3 -
			\cfrac{a_2 x^6}{ 1 + b_2 x^3 - \ddots}}}
\intertext{where\quad
	$a_n = (3n)(3n+1)^2 (3n+2)^2 (3n+3)$,\quad
	$b_n = 2(3n+3) \left( (3n+3)^2 + 1 \right)$.}
	\int_0^\infty \cm(u) \, e^{-u/x} \, du &=
		\frac{x}{1+2x^3+
			\cfrac{a_1 x^6}{ 1 + b_1 x^3 -
			\cfrac{a_2 x^6}{ 1 + b_2 x^3 - \ddots}}}
\intertext{where\quad
	$a_n = (3n-2)^2(3n-1)^2 (3n)^2$,\quad
	$b_n = (3n-1) (3n)^2 + (3n+1)^2 (3n+2) $.}
	\int_0^\infty \sm(u) \, \cm(u) \, e^{-u/x} \, du &=
		\frac{x^2}{1+12x^3+
			\cfrac{a_1 x^6}{ 1 + b_1 x^3 -
			\cfrac{a_2 x^6}{ 1 + b_2 x^3 - \ddots}}}
\intertext{where\quad
	$a_n = (3n-1)^2(3n)^2 (3n+1)^2$,\quad 
	$b_n = (3n) (3n+1)^2 + (3n+2)^2 (3n+3) $.}
	\int_0^\infty \sm^2(u) \, \cm(u) \, e^{-u/x} \, du &=
		\frac{2x^3}{1+40x^3+
			\cfrac{a_1 x^6}{ 1 + b_1 x^3 -
			\cfrac{a_2 x^6}{ 1 + b_2 x^3 - \ddots}}}
\intertext{where\quad
	$a_n = (3n)^2(3n+1)^2 (3n+2)^2$,\quad
	$b_n = (3n+1) (3n+2)^2 + (3n+3)^2 (3n+4) $.}
\end{align*}
\end{minipage}}
\caption{\label{sixJ-fig}
The six basic $J$-fractions relative to Dixonian functions.}
\end{figure}

Consider the list of positive integers in ascending order with each integer
listed twice:
\[
	1,\, 1,\, 2,\, 2,\, 3,\, 3,\, 4,\, 4,\, 5,\, 5,\, 6,\, 6,\, \dots
\]
There are essentially six ways to break this up into ascending 6-tuples,
allowing for the possibility of missing leading entries in the first
6-tuple.  The $a_n$ terms in the $J$-fraction for
the Laplace transform of $\sm(u)$, $\sm^2(u)$,  $\sm^3(u)$
are seen to correspond to the partitions:
\[
\begin{array}{ll}
\sm(u)~:& (1,2,2,3,3,4),\, (4,5,5,6,6,7),\, (7,8,8,9,9,10),\, \dots\,,
\\
\sm^2(u)~:& (2,3,3,4,4,5),\, (5,6,6,7,7,8),\, (8,9,9,10,10,11),\, \dots\,,
\\
\sm^3(u)~:& (3,4,4,5,5,6),\, (6,7,7,8,8,9),\, (9,10,10,11,11,12),\, \dots\,,
\end{array}
\]
%
%
These are the three possible ``odd'' partitions, odd in the sense that the
first integer in each 6-tuple appears exactly once. The three remaining
continued fractions, those of $\cm$, $\cm\cdot \sm$, $\cm\cdot \sm^2$,
 are associated in the same way with the three possible
even partitions.

\subsection{$S$-fractions for Dixon functions}
Starting with the original recurrences of the previous section
(Equation (\ref{Dixon:rec:first}), 
we can use the relation between $S_n$ and $D_{n-1}$
to eliminate either the letter $S$ or the
letter $D$.
These recurrences reduce to $S$-fraction recurrences, which
we tabulate here:
\[\begin{array}{lllllll}
\ds	\frac{S_n}{C_{n-2}} &=&\ds
		\cfrac{n(n-1)x^2}{ 1 + n(n+1)x^2 \frac{C_{n+1}}{S_n}}\,, &&
\ds	\frac{C_n}{S_{n-1}} &=&\ds
		\cfrac{n x}{ 1 + (n+1)x \frac{S_{n+2}}{C_n}} \\
\ds	\frac{C_n}{D_{n-2}} &=&\ds
		\cfrac{ n(n-1)x^2 }{ 1 + (n+1)(n+2)x^2 \frac{D_{n+1}}{C_n}}\,, &&
\ds	\frac{D_n}{C_{n-1}} &=&\ds
		\cfrac{n x}{ 1 + (n+2)x \frac{C_{n+2}}{D_n}}\,.
\end{array}
\]

We need initial conditions to generate $S$-fractions, and the six candidate
starting points give just three $S$-fraction initial conditions:
\[\begin{array}{lllll}
	S_1 &=& x^2 - 2x^2 C_2
		& = &\ds\cfrac{x^2}{1 + 2x \frac{C_2}{S_1}} \\
	C_0 &=& x - x S_2
		& = &\ds\cfrac{x}{1 + x \frac{S_2}{C_0}} \\
	C_1 &=& x^2 - 2x S_3
		& = &\ds\cfrac{x}{1 + 2x \frac{S_3}{C_1}}\,.
\end{array}
\]
The three that fail to give good initial conditions for an $S$-fraction
are as follows:
\[
\begin{array}{llll}
	S_2 &= 2x^3 - 2x^3 S_2 - 6x^2 C_3
		&& = \cfrac{2x^3}{1 + 2x^3 + 6x^2 \frac{C_3}{S_2}} \\
	S_3 &= 6x^4 - 12x^3 S_3 - 12x^2 C_4
		&& = \cfrac{6x^4}{1 + 12x^3 + 12x^2 \frac{C_4}{S_3}} \\
	C_2 &= 2x^3 - 4x^3 C_2 - 3x^2 S_4
		&& = \cfrac{2x^3}{1 + 4x^3 + 3x^2 \frac{S_4}{C_2}} \,.\\
\end{array}
\]

On iterating the first three relations, we obtain three
$S$-fraction expansions:

\begin{theorem}[Conrad~\cite{Conrad02}] The formal Laplace transform
of the functions $\sm$, $\cm$ and $\sm\cdot \cm$ have explicit 
$S$-fraction expansions with cubic numerators,
\begin{align*}
	\int_0^\infty \sm u \, e^{-u/x} \, du &=
	\frac{x^2}{1 + 
		\cfrac{a_1 x^3}{1 +
		\cfrac{a_2 x^3}{1 + \ddots }}} \\
\intertext{where for $r\geq 1$:\quad
	$a_{2r-1} = (3r-2)(3r-1)^2$, \quad
	$a_{2r} = (3r)^2 (3r+1)$;}
	\int_0^\infty \cm u \, e^{-u/x} \, du &=
	\cfrac{x}{1 + 
		\cfrac{a_1 x^3}{1 +
		\cfrac{a_2 x^3}{1 + \ddots }}} \\
\intertext{where for $r\geq 1$:\quad
	$a_{2r-1} = (3r-2)^2 (3r-1)$,\quad
	$a_{2r} = (3r-1) (3r)^2$;}
	\int_0^\infty \sm(u) \, \cm(u) \, e^{-u/x} \, du &=
	\cfrac{x^2}{1 +
		\cfrac{a_1 x^3}{1 +
		\cfrac{a_2 x^3}{1 + \ddots }}} \\
\intertext{where for  $r\geq 1$:\quad
	$a_{2r-1} = (3r-1)^2 (3r)$,\quad
	$a_{2r} = (3r) (3r+1)^2$.}
\text{.}
\end{align*}
\end{theorem}

Consider now the set of positive integers written in increasing order with
each integer written twice:
\[
	1,\, 1,\, 2,\, 2,\, 3,\, 3,\, 4,\, 4,\, 5,\, 5,\, 6,\, 6,\, \dots
\]
There are three ways to divide this list into ascending triples if we
permit ourselves to discard leading ones, corresponding to the
numerators coefficients $a_n$ of the $S$-fractions associated to 
$\cm$, $\sm$, and $\sm\cdot\cm$:
\[
\begin{array}{ll}
\cm(u)~: & (1,1,2),\, (2,3,3),\, (4,4,5),\, (5,6,6),\, \dots\,,
\\
\sm(u)~: & (1,2,2),\, (3,3,4),\, (4,5,5),\, (6,6,7),\, \dots\,,
\\
\sm(u)\cdot\cm(u)~: & (2,2,3),\, (3,4,4),\, (5,5,6),\, (6,7,7), \,
\dots\,.
\end{array}
\]
This correspondence
suggests that these three continued fractions form a complete set.

\smallskip

\begin{Note} \emph{The cubic $X^3+Y^3-3\alpha XY=1$.}
Conrad in his dissertation~\cite{Conrad02}, follows Dixon
and examines the larger class of functions
corresponding to the cubic
\[
X^3+Y^3-3\alpha XY=1,
\]
for arbitrary $\alpha$.  These do give rise to continued fractions, but
ones that are non-standard:
they are not of the $S$ or $J$ type as they involve 
denominators that are linear (in~$x$) and numerators that are cubic.

\end{Note}

An alternative derivation of the $J$-fraction expansions, based on a
direct use of the differential system, will be given in
Subsection~\ref{confrac-sec}, when we discuss a method of Andr\'e
(Note~\ref{andre-note}). 

\section{First combinatorial model: Balls games} \label{balls-sec}

Dixonian elliptic functions, as we'll soon see, serve to
describe the evolution
of a simple urn model with balls of two
colours. In this perspective, the Taylor coefficients of $\sm,\cm$
count certain combinatorial objects that are urn histories,
or equivalently, weighted knight's walks in the discrete plane---this provides our \emph{first
combinatorial interpretation}. Going from the discrete to the
continuous, we furthermore find that Dixonian functions quantify the 
extreme behaviour 
of a classic continuous-time branching process, the Yule process. 
Finally, we show that the composition of 
the system at any instant, whether in the discrete or the continuous,
case, can be fully worked out and is once more expressible in terms of Dixonian functions.

\subsection{Urn models.} \label{Polya-sec}
Balls games have been of interest to probabilists since the dawn of
time. For instance in his \emph{Th\'eorie analytique des
probabilit\'es} (first published in~1812), Laplace writes
``\emph{Une urne~$A$ renfermant un tr\`es grand nombre~$n$ de boules
blanches et noires; \`a chaque tirage, on en extrait une que l'on
remplace par une boule noire; on demande la probabilit\'e
qu'apr\`es~$r$ tirages, le nombre des boules blanches sera~$x$}.''
Such games were later systematically studied by P\'olya. 
What is directly relevant to us
here is the following version known in the standard probability
literature as the \emph{P\'olya urn model}
(also P\'olya-Eggenberger):
\begin{itemize}
\item[]\em
{\bf P\'olya urn model.}
An urn is given that contains black and white balls. 
At each epoch, a ball in the urn
is chosen at random (but \emph{not} removed). If it is black, then $\alpha$ black and $\beta$
white balls are placed into the urn; else it is white and $\gamma$ black and
$\delta$ white balls are placed into the urn.
\end{itemize}
The model is fully described by the ``placement matrix'',
\[
\cal M=\left(\begin{array}{cc}\alpha&\beta\\ \gamma&\delta\end{array}\right).
\]
The most frequently encountered models are \emph{balanced}, meaning
that $\alpha+\beta=\gamma+\delta$, and
negative entries in a
matrix~$M$ are interpreted as
subtraction (rather than addition) of balls. 
For instance, Laplace's original problem corresponds to
{\small\renewcommand{\arraycolsep}{2pt}\renewcommand{\arraystretch}{0.7}
 $\left(\begin{array}{cc}0&0\\1&-1\end{array}\right)$},
which is none other than the coupon collector's problem in modern terminology.
What is sought in various areas of science
is some characterization, exact or asymptotic, of
the composition of the urn at epoch~$n$, given fixed initial
conditions. The elementary introduction by
Johnson and Kotz~\cite{JoKo77} mentions applications to
sampling statistics, learning processes, decision theory, and
genetics. Recently, P\'olya urn models have been found to be of
interest in the analysis of several algorithms and data structures of
computer science; see especially Mahmoud's survey~\cite{Mahmoud03}.

\smallskip

The main character of this section is the special urn defined by the
matrix
\[
\M_{12}=\left(\begin{array}{cc}-1&2\\2&-1\end{array}\right).
\]
This can be visualized as a game with balls of either black (`$\b$') or
white (`$\w$') colour. If a ball is chosen, it is removed [the $-1$
entry in the matrix] from the
urn and replaced by two balls of the opposite colour
according to the rule
\begin{equation}\label{rule12}
\b~\longrightarrow~\w\w,\qquad  
\w~\longrightarrow~\b\b.
\end{equation}

A \emph{history} of length~$n$
(see Fran{\c c}on's work~\cite{Francon78} for this terminology and
similar ideas) 
is, loosely speaking, any description of a
legal sequence of $n$ moves of the
P\'olya urn. Let conventionally the urn be initialized at time~0 with
one black ball~({\sf x}). A history (of length~$n$)
is obtained by starting with the
one-letter word $\b$ at time~0 
and successively applying ($n$ times) the rules
of~(\ref{rule12}). For instance,
\[
\underline{\b} \longrightarrow \w\underline{\w} \longrightarrow 
\w\underline{\b}\b \longrightarrow
\underline{\w}\w\w\b \longrightarrow 
\b\underline{\b}\w\w\b\longrightarrow \b\w\w\w\w\b,
\]
is the complete description of a history of length~$5$.
(The replaced letters have been underlined for readability.)
We let $H_{n,k}$ be the number of histories that start with
an~$\b$ and, after~$n$ actions, result in a word having $k$
occurrences of~$\w$ (hence $n+1-k$ occurrences of~$\b$).
Clearly the total number of histories of length~$n$
satisfies
\[
H_n:=\sum_k H_{n,k}=n!,
\]
since the number of choices is $1,2,3,\ldots$ at times
$1,2,3\ldots$\,.
Here is a small table of all histories of length
$\le3$:
\begin{small}
\[
\begin{array}{lll}
n=0~:&& \x \\
n=1~:&& \x \lora \y\y \\
n=2~:&& \x \lora \y\y \lora \x\x\y \\
&& \x \lora \y\y \lora \y\x\x \\
n=3:~&& \x \lora \y\y \lora \x\x\y \lora \y\y\x\y\\
&& \x \lora \y\y \lora \x\x\y \lora \x\y\y\y \\
&& \x \lora \y\y \lora \x\x\y \lora \x\x\x\x \\
&& \x \lora \y\y \lora \y\x\x \lora \x\x\x\x \\
&& \x \lora \y\y \lora \y\x\x \lora \y\y\y\x \\
&& \x \lora \y\y \lora \y\x\x \lora \y\x\y\y
\end{array}
\]
\end{small}%
The sequence $(H_{n,0})$ starts as $1,0,0,2$ for
$n=0,1,2,3$ and it is of interest to characterize these combinatorial
numbers. 

\subsection{Urns and Dixonian functions.}\label{urns-sec}
Let us come back to Dixonian functions. 
Consider for notational convenience 
the (autonomous, nonlinear) ordinary differential system
\begin{equation}\label{sysu}
\Sigma: \qquad \frac{dx}{dt}=y^2,\quad \frac{dy}{dt}=x^2,\quad
\hbox{with}\quad x(0)=x_0,\quad y(0)=y_0,
\end{equation}
which is the signless version of~(\ref{conrad1}).
In this subsection, we only need the specialization
$x(0)=0$, $y(0)=1$, but we will make use of the general case~(\ref{sysu}) in
Subsection~\ref{dyn-sec} below.

The pair $\langle x(t),y(t)\rangle$
with initial conditions $x(0)=0$, $y(0)=1$ parametrizes the ``Fermat
hyperbola'', 
\[
y^3-x^3=1,
\]
which is plainly obtained from Fermat's ``circle'' $\F_3$ by a
vertical symmetry.
It is natural to denote this pair
of functions by $\langle\smh(t),\cmh(t)\rangle$, 
these functions $\smh,\cmh$ being 
trivial variants of $\sm,\cm$, where the alternation of signs has been
suppressed:
\begin{equation}\label{scmhdef}
\smh(z)=-\sm(-z),\qquad \cmh(z)=\cm(-z).
\end{equation}


A first combinatorial approach to Dixonian functions can be
developed in a straightforward manner by simply looking at the 
basic algebraic relations induced by the system~$\Sigma$ of~(\ref{sysu}).
To this purpose,
we define a linear transformation
$\delta$ acting on the vector space
$\C[x,y]$ of polynomials in two formal variables $x,y$ that is specified by
the rules,
\begin{equation}\label{rules}
\delta[x]=y^2, \qquad \delta[y]=x^2,\qquad
\delta[u\cdot v]=\delta[u]\cdot
v + u\cdot \delta[v],
\end{equation}
$u,v$ being arbitrary elements of $\C[x,y]$.
A purely \emph{mechanical} way to visualize the operation of~$\delta$
comes from regarding $\delta$ as a rewriting
system
\begin{equation}\label{rew}
x\tode yy, \qquad y \tode xx,
\end{equation}
according to the first two rules of~(\ref{rules}).
The third rule means that $\delta$ is a derivation, and it can be read
algorithmically as follows: when $\delta$ is to be applied to a
monomial~$w$ of total degree~$d$ in $x,y$, first arrange $d$ copies 
of $w$, where in each copy one instance of a variable is marked
(underlined), then apply the rewrite rule~(\ref{rew}) once in 
each case to the marked variable, and finally collect the results. For
instance, 
\[
xyy \mapsto \underline{x}yy, x\underline{y}y, xy\underline{y}
\tode (yy)yy, x(xx)y, xy(xx), \quad \hbox{so  that}\quad
\delta[xy^2]=y^4+2x^3y.
\]

Two facts should now be clear from the description of~$\delta$
and the definition of histories.
\begin{itemize}
\item[$(i)$] Combinatorially, the $n$th iterate $\delta^n[x^ay^b]$
describes the 
collection of all the possible histories at time~$n$ 
of the P\'olya urn with matrix~$\M_{12}$,
when the initial configuration of the urn has $a$  balls of the first
type ($\x$) and $b$ balls of the second type ($\y$). 
This is exactly the meaning of the replacement
rule~(\ref{rew}). In particular, the
coefficient\footnote{%
	If $f=\sum_{m,n}f_{m,n}x^my^n$, then the notation
$[x^my^n]f$ is used to represent coefficient extraction:
$[x^my^n]f\equiv f_{m,n}$.
}, 
\[
H_{n,k}^{(a,b)} = [x^k y^{\ell}]\, \delta^n[x^ay^b],\qquad
k+\ell=n+a+b
\]
is the number of histories of a P\'olya urn
that lead from the initial state $\x^ay^b$ to the final state~$\x^k\y^\ell$.
\item[$(ii)$] 
Algebraically, the operator~$\delta$ does nothing but describe the
``logical consequences'' 
of the differential system $\Sigma$. In effect, the first two rules
of~(\ref{rew}) mimic the effect of a derivation 
applied to terms containing $x=x(t)$ and $y=y(t)$, ``knowing'' that
$x'=y^2$, $y'=x^2$. Accordingly, the quantity
$\delta^n[x^ay^b]$ represents an $n$th derivative,
\[
\delta^n[x^ay^b]
=\frac{d^n}{dt^n} x(t)^a y(t)^b\quad\hbox{expressed in}\quad x(t),y(t),
\]
where $x(t),y(t)$ solve the differential system $x'=y^2$, $y'=x^2$.
\end{itemize}
In summary, we have a principle\footnote{%
  In this paper, we have chosen to develop a calculus geared to 
  Dixonian functions. Our approach is in fact appreciably more general
  and it can be applied to any balanced urn model,
 as considered in~\cite{FlGaPe05}.}
 whose informal statement
is as follows.
\begin{itemize}
\item[]
{\bf Equivalence principle.}
\emph{The algebra of a nonlinear autonomous system that is monomial
and homogeneous, like~$\Sigma$
in~\eqref{sysu},
is isomorphic to the combinatorics of an associated P\'olya urn.}
\end{itemize}
This principle can be used in either direction. 
For us it makes it possible to analyse the P\'olya urn $\M_{12}$
in terms of the functions $\sm,\cm$ that have already been made explicit in
Proposition~\ref{prop-hyper}. In so doing, we rederive and expand
upon an analysis given in a recent paper on analytic and probabilistic
aspects of urn models~\cite{FlGaPe05}.
Here is one  
of the many consequences of this equivalence.

\begin{theorem}[First combinatorial interpretation] \label{urn0-prop}
The exponential generating function of histories 
of the urn with matrix $\M_{12}$ that start with one ball
and terminate with balls that are \emph{all} of the other colour is
\[
\sum_{n\ge0} H_{n,0}\frac{z^n}{n!}
=\smh(z)=\frac{\sm(z)}{\cm(z)}=-\sm(-z)
.
\]
The exponential generating function of histories 
that start with one ball and terminate with balls that are all of
the original colour is
\[
\sum_{n\ge0} H_{n,n+1}\frac{z^n}{n!}=\cmh(z)=\frac{1}{\cm(z)}=\cm(-z).
\]
\end{theorem}
\begin{proof}
This is nothing but Taylor's formula. Indeed, for the first equation,
we have by the combinatorial interpretation of the differential
system:
\[
H_{n,0}=\left.\delta^n [x]\right|_{x=0,y=1}.\]
But from the algebraic interpretation,
\[
\left.\delta^n[x]\right|_{x=0,y=1} = 
\left.\frac{d^n}{dt^n} \smh(t)\right|_{t=0}.\]
The result then follows.
\end{proof}

Theorem~\ref{urn0-prop} thus provides a first combinatorial model of
Dixonian functions in terms of urn histories.

\begin{Note} \emph{On Dumont.} 
Though ideas come from different sources, there is a striking parallel
between what we have just presented and some of Dumont's researches in the
1980's  and 1990's. Dumont gives an 
elegant presentation of Chen grammars in~\cite{Dumont96}. There, he
considers chains of general substitution rules on words: such chains
are  partial differential operators in disguise.
For instance, our $\delta$
operator is nothing but
\[
\delta = x^2\frac{\partial}{\partial y} + y^2\frac{\partial}{\partial
x}.\]
Dumont has shown in~\cite{Dumont96} that operators of this type
can be used to
approach a variety of questions like rises in permutations, Bell
polynomials, increasing trees, parking functions, and P\'olya
grammars. The interest  of such investigations 
is also reinforced by consideration of the trivariate operator
\[
yz\frac{\partial}{\partial x}+zx\frac{\partial}{\partial y}
+xy\frac{\partial}{\partial z},
\]
itself related to the differential system
\[
x'=yz,\quad y'=zx,\quad z'=xy,
\]
and to elliptic functions of the Jacobian type, which had been formerly
researched by Dumont and Schett; see~\cite{Dumont79,Dumont81,Schett76}
and Subsection~\ref{sixlat-sec} below for further comments on the
operator point of view.
\end{Note}

\begin{figure}
\setlength{\unitlength}{1.truemm}\thicklines
\begin{picture}(30,30)
\put(10,10){\circle*{2}}
\put(-8,10){$P=(p,q)$}
\put(10,10){\line(0,1){20}}
\put(10,20){\circle*{1}}
\put(10,30){\circle*{1}}
\put(10,30){\vector(-1,0){10}}
\put(0,30){\circle*{2}}
\put(10,10){\line(1,0){20}}
\put(20,10){\circle*{1}}
\put(30,10){\circle*{1}}
\put(30,10){\vector(0,-1){10}}
\put(30,0){\circle*{2}}
\end{picture}
\hspace*{1.5truecm}
\begin{picture}(30,30)
\put(10,10){\circle*{2}}
\put(-8,10){$P=(p,q)$}
\put(10,10){\line(0,1){20}}
\put(10,20){\circle*{1}}
\put(12,25){multiplicity $p$}
\put(10,30){\circle*{1}}
\put(10,30){\vector(-1,0){10}}
\put(0,30){\circle*{2}}
\put(10,10){\line(1,0){20}}
\put(20,10){\circle*{1}}
\put(25,12){multiplicity $q$}
\put(30,10){\circle*{1}}
\put(30,10){\vector(0,-1){10}}
\put(30,0){\circle*{2}}
\end{picture}
\caption{\label{moves-fig} The Bousquet-M\'elou-Petkov{\v s}ek moves
(left) and their version with multiplicities (right) for the $\M_{12}$
urn model.}
\end{figure}
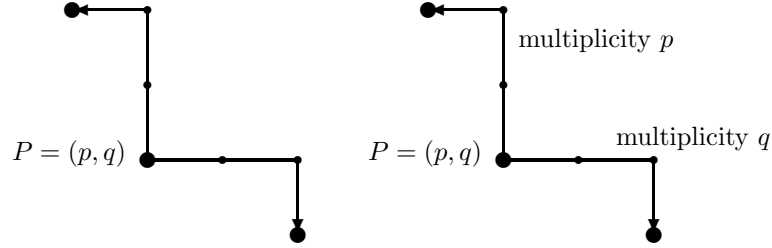

\begin{Note} \emph{Knight's walks of Bousquet-M\'elou and Petkov{\v s}ek.} 
Theorem~\ref{urn0-prop} is somehow related to other interesting
combinatorial objects (Figure~\ref{moves-fig}). 
Say we only look at the possible evolutions
of the urn, disregarding which particular ball is chosen.
An urn with $p$ black balls and $q$ white balls may be represented in
the cartesian plane by the point $P$ with coordinates~$(p,q)$. 
A sequence of moves then defines a polygonal line, $P_0P_1\cdots P_n$
with $P_0=(1,0)$, each move $\overrightarrow{P_{j}P}_{j+1}$ being of type either
$\beta=(-1,+2)$ [i.e., two steps East, one step South] when a black
ball is picked up
or $\omega=(2,-1)$ [i.e., two steps North, one step West] otherwise.
This defines a random walk on the first quadrant that makes use
of two types of knight moves on the chessboard. The enumeration of
these walks is a nontrivial combinatorial problem that has been solved
recently by Bousquet-M\'elou and Petkov{\v s}ek in~\cite{BoPe03}. 
They show for instance that the ordinary generating function of walks
that start at $(1,0)$ and end on the horizontal axis is
\[
G(x)=\sum_{i\ge0} (-1)^i \left(\xi^{\langle i\rangle}(x)\xi^{\langle
i+1\rangle}(x)\right)^2,
\]
where $\xi$ is a branch of the (genus 0) cubic $x\xi-x^3-\xi^3=0$:
\[
\xi(x)=x^2\sum_{m\ge0} \binom{3m}{m}\frac{x^{3m}}{2m+1}.
\]
There, $\xi^{\langle i\rangle}=\xi\circ\cdots\circ\xi$ is the $i$th
iterate of~$\xi$. These walks have a merit~\cite{BoPe03}: they provide 
an extremely simple example of a linear recurrence with constant
coefficients whose generating function is highly transcendental (in fact,
not even holonomic).

In order to obtain a
complete history from a knight's walk, one has to add some
supplementary information, 
namely, which ball is chosen at each stage. This corresponds to
introducing multiplicative weights. The rule is then as follows. For a walk
with a $\beta$-step that starts at point $(p,q)$
(a black ball is chosen) the weight is $p$; dually, for an
$\omega$-step the weight 
is~$q$. 
\end{Note}

In summary, at this stage, we have available three variants
of the interpretation of the Taylor coefficients of the Dixonian functions
$\sm,\cm$ provided by Theorem~\ref{urn0-prop}: $(i)$~the enumeration of urn histories relative to the $\M_{12}$
urn; $(ii)$~the iterates of the special operator
$\delta=x^2\partial_y+y^2\partial_x$ in the style of Dumont; 
$(iii)$~multiplicatively weighted knight walks of the type introduced by
Bousquet-M\'elou and Petkov{\v s}ek.

\begin{Note} \emph{A probabilistic consequence: extreme large deviations.}
From a probabilistic standpoint, 
the P\'olya urn model is a discrete time Markov chain with a
denumerable set of states embedded in~$\Z\times\Z$.
The number of black balls at time~$n$
then becomes a random variable, $X_n$.
Theorem~\ref{urn0-prop} quantifies the probability ($\Pr$)
 of extreme large deviations of~$X_n$ 
as
\[
\Pr(X_n=0)=\frac{H_{n,0}}{n!}=[z^n]\smh(z)=[z^n]-\sm(-z).
\]
Then by an easy analysis of singularities,
one finds that this quantity decreases exponentially fast:
\[
\Pr(X_{3\nu+1}=0) \sim c\left(\frac{\pi_3}{3}\right)^{-3\nu-1},\qquad
\hbox{where}\quad \frac{\pi_3}{3}=\frac{\sqrt{3}}
{6\pi}\Gamma\left(\frac13\right)^3,
\] 
in accordance with~(\ref{pi3def}).
For instance, we have to 10D:
\[
\left|\frac{[z^{28}]\sm(z)}{[z^{31}]\sm(z)}\right|^{1/3}\doteq
1.76663\,87{\it 502}\cdots;
\qquad
\frac{\pi_3}{3}\doteq
1.76663\,87490\cdots\,.
\]
One can refer to the calculation of~\cite{FlGaPe05}, but is is just as easy
to note that Proposition~\ref{prop-hyper} provides directly the
dominant positive 
singularity $\rho$ of~$\smh(z)$ as a special value of the
fundamental Abelian integral. 
\end{Note}

\subsection{A continuous-time branching process.} \label{Yule-sec}
\begin{figure}[b]
\begin{center}
\Img{5}{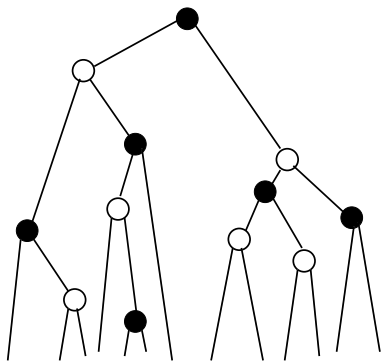}
\end{center}
\caption{\label{yule-fig}
A truncated view of a tree of particles provided by the Yule process
(the ordinate of each particle represents the time at which it
splits into two particles of the other type).}
\end{figure}

The Dixonian functions also make it possible to answer questions
concerning chain reactions in a certain form of
particle physics. You have two types of particles, say,
foatons and viennons. Any particle lives an amount of time~$T$ that is
an exponentially distributed random variable
(i.e., $\Pr(T\ge t)=e^{-t}$), this
independently of the other particles; then it
disintegrates into two particles of the other type. Thus a foaton
gives rise to two viennons and a viennon gives rise to two foatons
(see Figure~\ref{yule-fig}).
What is the composition of the system at some time~$t\ge0$, assuming
one starts with one foaton? 

The resemblance with the P\'olya urn should be clear (with a foaton
being a black~$\x$ and a viennon a white~$\y$).
First, considering only the total population of particles, 
one gets what is perhaps the simplest continuous-time branching process,
classically known as the \emph{Yule process}. 
(For useful information regarding this process, see for instance the recent
articles by Chauvin,
Rouault, and collaborators~\cite{ChKlMaRo04,ChRo04}.) Let $S_k(t)$ be
the probability 
that the total population at time $t$ is of size $k$. Introduce the
bivariate generating function,
\[
\Xi(t;w)=\sum_{k=1}^\infty S_k(t) w^k.
\]
Examination of what happens between times~0 and $dt$ leads to
the recurrence
\begin{equation}\label{forgot}
S_k(t+dt)=(1-dt)S_k(t)+dt \sum_{i+j=k}S_i(t)S_j(t),
\end{equation}
which is none other than the usual ``backwards equation'' of Markov
processes. Then, one has
\[
S_k'(t)+S_k(t)=\sum_{i+j=k}S_i(t)S_j(t).\]
This induces a nonlinear equation satisfied by $\Xi$,
namely,
\[
\Xi'(t;w)+\Xi(t;w)=\Xi(t;w)^2, \qquad \Xi(0,w)=w,
\]
where derivatives are implicitly taken with respect to the time
parameter~$t$. The solution of this ordinary differential equation is
easily found by separation of variables,
\[
\Xi(t;w)=\frac{we^{-t}}{1-w(1-e^{-t})},
\]
which yields 
\begin{equation}\label{yulesize}
S_k(t)=e^{-t}\left(1-e^{-t}\right)^{k-1}, \qquad k\ge1.
\end{equation}
In summary, the size of the
population at time~$t$ obeys a geometric law of
parameter $(1-e^{-t})$, with expectation  $e^t$. 
(The previous calculations are of course extremely classical: they
are disposed of in just six lines in Athreya and Ney's
treatise~\cite[p.~109]{AtNe72}.) 

The result of~(\ref{yulesize}) shows that any calculation under the
discrete urn model can be automatically transferred to the continuous branching
process. 
Precisely, let $\cal H$ be the set of all histories of the P\'olya urn with
matrix ${\cal M}_{12}$.  Let $\cal K \subseteq \cal H$ be a subset of
$\cal H$ and let $K(z)$ be the exponential generating function of
$\cal K$. Then, by virtue of~(\ref{yulesize}),
the probability that, at time~$t$, the Yule process
has evolved according to a history that lies in $\cal K$ is given by
\begin{equation}\label{discrcont}
e^{-t}K\left(1-e^{-t}\right).
\end{equation}
This has an immediate consequence.


\begin{proposition} \label{Yule-prop}
Consider the Yule process with two types of particles.
The probabilities that particles are all of the second type 
at time~$t$ are 
\[
X(t)=e^{-t}\smh(1-e^{-t}),\qquad
Y(t)=e^{-t}\cmh(1-e^{-t}),
\]
depending on whether the system at time~0 is
initialized with one particle of the
first type ($X$) or of the second type ($Y$).
\end{proposition}

Like its discrete-time counterpart,
this proposition can  be used to quantify extreme deviations:
the probability that all particles be of the second type at time~$t$ is
asymptotic to
\begin{equation}\label{single}
e^{-t}\smh(1),\qquad \smh(1)\doteq 1.20541\,51514.
\end{equation}
At time~$t$, the Yule system has expected size $e^t$. On the other hand, 
for the discrete process, the
probability at a large discrete time that balls 
or particles be all of one colour 
is exponentially small. These two facts might lead to
expect that, in continuous time, the probability of all particles to be
of the same colour is doubly (and not singly, as in
Equation~(\ref{single})) exponentially small. The apparent paradox is
resolved by 
observing that the extreme large deviation regime is driven by the
very few cases where the system consists of $O(1)$ particles only, an event
whose probability is exponentially small (but \emph{not} a double exponential). 

\begin{Note} \emph{An alternative direct derivation.}
The memoryless nature of the process
implies, like in~(\ref{forgot}), the differential system:
\[
X(t)+X'(t)=Y(t)^2,\qquad Y(t)+Y'(t)=X(t)^2,
\]
with initial conditions $X(0)=0$, $Y(0)=1$. This system  closely
resembles some of our earlier equations. By simple algebra [first multiply by
$e^t$, then apply the change of variables $t\mapsto e^{-t}$],
the solution in terms of Dixonian functions results.
\end{Note}

\subsection{Dynamics of the P\'olya and Yule processes.}\label{dyn-sec}
The previous sections have quantified the extreme behaviour of the
processes---what is the probability that balls be all of one colour?
In fact, a slight modification of previous arguments give 
complete access to the
composition of the system at any instant.

Consider once more the differential system 
\begin{equation}\label{sysubis}
\frac{dx}{dt}=y^2,\quad \frac{dy}{dt}=x^2,\quad
\hbox{with}\quad x(0)=x_0,\quad y(0)=y_0
\end{equation}
(this repeats Equation~(\ref{sysu})). 
The initial conditions are now treated as free parameters
or formal variables. This system can be solved
exactly by means of previously exposed techniques, and,
by virtue of the equivalence principle,
its  solution describes combinatorially the composition of the urn, not
just the extremal configurations (that correspond to $x_0=0$, $y_0=1$).

First, the solution to the system. Following the chain of~(\ref{sol1})
and~(\ref{sol2}), one finds
\[
x'(t)^{3/2}=x(t)^3+\Delta^3, \qquad \Delta^3:=y_0^3-x_0^3,\]
where the initial conditions have been taken into account. 
This last equation can be put under the form
\[
\frac{x'}{(x^3+\Delta^3)^{2/3}}=1,\]
which integrates to give
\begin{equation}\label{zz}
\Delta t=\int_{x_0/\Delta}^{x/\Delta} \frac{dw}{(1+w^3)^{2/3}}.
\end{equation}
The function~$\smh(t)$ is also defined by an inversion of an integral,
\[
t=I(\smh(t)), \qquad \hbox{where} \quad
I(y):=\int_0^y \frac{dw}{(1+w^3)^{2/3}},
\]
as is easily verified by the technique of Section~\ref{hyper-sec}.
Then Equation~(\ref{zz}) provides
\begin{equation}\label{zz2}
x(t)=\Delta \smh\left(\Delta t +
I\left(\frac{x_0}{\Delta}\right)\right),
\end{equation}
which constitutes our main equation.

Second, the interpretation of~$\delta$ in Section~\ref{urns-sec}
implies that
$\delta^n[x]$ expresses the composition of the urn at the
$n$th stage of its operation when it has been initialized with
one ball of the first type ($\x$). Thus the quantity 
\begin{equation}\label{zz3}
F(t,x_0)=\sum_{n\ge0} \frac{t^n}{n!} \left(\delta^n[x]\right)_{x\mapsto x_0,\,y\mapsto1}
\end{equation}
is none other than the $x(t)$ solution to the system~(\ref{sysubis})
initialized with $x(0)=x_0$ and $y(0)=1$.
On the other hand, $\delta^n[x]$ represents
the $n$th derivative of $x\equiv x(t)$ expressed as
a function of~$x$ and~$y$. 

There results from Equations~(\ref{zz2}), (\ref{zz3}) and the
accompanying remarks an expression
for the bivariate generating function of 
histories with~$t$ marking length and~$x_0$ marking the number of
balls of the first kind. Switching to more orthodox
notations ($x_0\mapsto x$), we state:

\begin{figure}
\begin{center}
\Img{5.5}{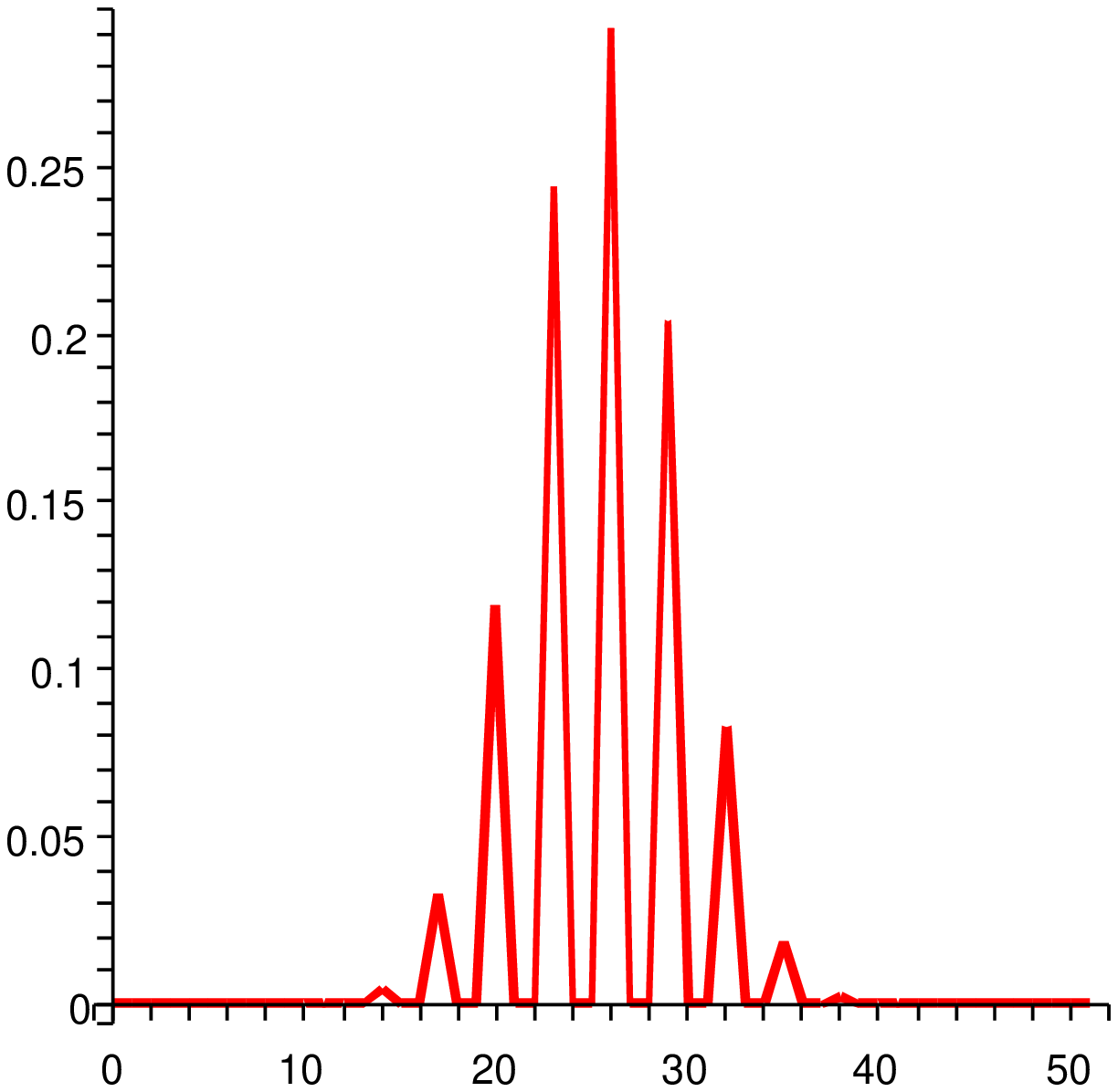} \Img{5.5}{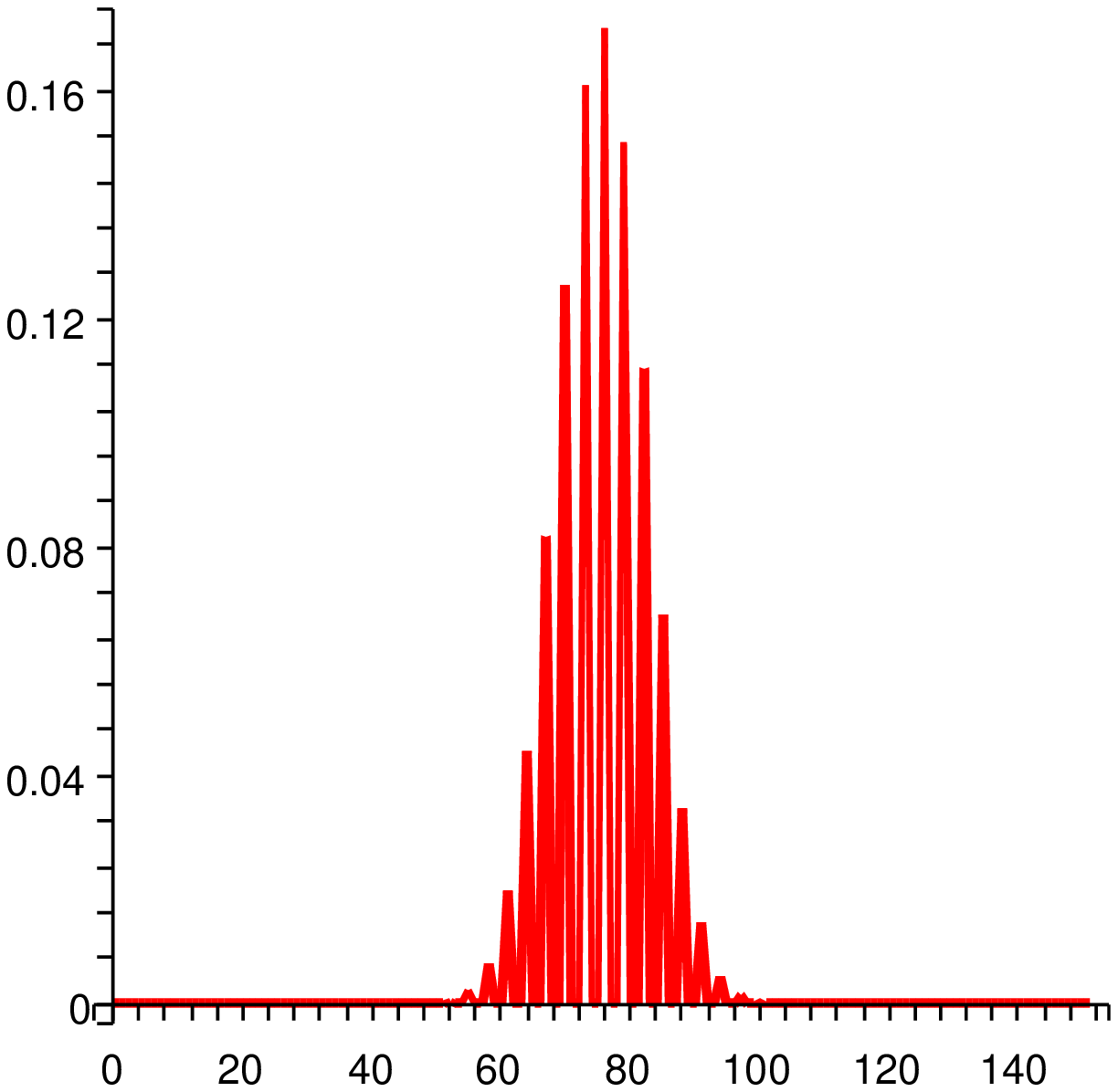}
\end{center}
\caption{\label{gauss-fig}
A plot of the coefficients of $\frac{1}{n!}\delta^n[x]$ for $n=50$ and
$n=150$ illustrating the asymptotically gaussian character of coefficients.
}
\end{figure}

\begin{proposition} \label{alltimes-prop}
The composition of the $\cal M_{12}$ urn at all
discrete instants  is described by
\[
\begin{array}{lll}
\ds \sum_{n,k\ge0} H_{n,k}x^k \frac{z^n}{n!}
&=&\ds (1-x^3)^{1/3}\smh\left((1-x^3)^{1/3}z+\int_0^{x(1-x^3)^{-1/3}}\!\!\!\!\!\!
\frac{ds}{(1+s^3)^{2/3}}\right),
\\
&=& \ds
\left(x+y^2\frac{z}{1!}+2x^2y\frac{z^2}{2!}+(4xy^3+2x^4)\frac{z^3}{3!}+\cdots
\right)_{y=1},
\end{array}
\]
where $H_{n,k}$ is the number of histories of length~$n$ of the urn
initialized with~$\x$ and terminating with $\x^k\y^{n+1-k}$.
\end{proposition}
This statement is the incarnation of Theorem~1
of~\cite{FlGaPe05} in the case of the $\cal M_{12}$ urn model. 
From it, one can for instance deduce by an application of the
Quasi-powers Theorem of analytic combinatorics~\cite{FlGaPe05,FlSe04,Hwang98b}:
\emph{The
distribution of the number of balls of the first kind,
equivalently the sequence of coefficients of~ the homogeneous
polynomial $\delta^n[x]$, is
asymptotically normal.} See Figure~\ref{gauss-fig}.

\smallskip

The transfer from discrete to continuous time afforded
by~(\ref{discrcont}) now permits us to deduce the composition of the
Yule process from the P\'olya urn via the transformation $K(z)\mapsto
e^{-t}K(1-e^{-t})$.

\begin{proposition} In the Yule process with two types
of particles, the probability generating function
of the number of particles of the first type $(\x$)
 at time~$t$ is 
\[
e^{-t}(1-x^3)^{1/3}\smh\left((1-x^3)^{1/3}(1-e^{-t})+\int_0^{x(1-x^3))^{-1/3}}
\!\!\!\!\!\!\frac{ds}{(1+s^3)^{2/3}}\right).
\]
\end{proposition}
By means of the continuity theorem for
characteristic functions, it can be verified that the distribution of
the number of particles of the first type
at large times is asymptotically exponential with mean $\sim e^t/2$.
(In this case, the distribution is essentially driven by the size of the system.)

\subsection{Complements regarding elliptic urn models.} \label{sixlat-sec} 
As already noted, the operator $\delta$ of the $\M_{12}$ urn
is a partial differential operator,
\[
\delta[f]=y^2\frac{\partial}{\partial x}f+x^2\frac{\partial}{\partial
y}f,
\]
which is linear and first order. 
In a recent study, Flajolet, Gabarr\'o and Pekari~\cite{FlGaPe05},
investigate the general class of 
$2\times2$ balanced urn schemes corresponding to matrices of the form
\begin{equation}\label{Mgen}
\M=\left(\begin{array}{cc}-a&s+a\\s-b&b\end{array}\right).
\end{equation}
Any such urn is modelled by a particular partial differential
operator,
\begin{equation}\label{gam}
\delta =x^{1-a}y^{s+a}\frac{\partial}{\partial x}+x^{s+b}
y^{1-b}\frac{\partial}{\partial
y},\qquad (a,b,s>0).
\end{equation}
\begin{figure}\small
\newcommand{\urn}[4]{\binom{#1~#2}{#3~#4}}
\begin{center}
\begin{tabular}{l|l|l|l|l}
\hline
& \hfil ($\cal M$) \hfil & \hfil  ($\Sigma$)\hfil &  
\hfil ($\delta$)\hfil  & \em Type
\\
\hline
$A$ & $\urn{-2}{3}{4}{-3}$ & $x'=x^{-1}y^3,~y'=x^4y^{-2}$ &
$x^{-1}y^3\der{x}+x^4y^{-2}\der{y}$ & $\wp$ \cite{FlGaPe05,PaPr98b};
\\
&&&& also Dixonian (\S\ref{jcwei-subsec})\\
$B$ &  $\urn{-1}{2}{3}{-2}$ & $x'=y^2,~y'=x^3y^{-1}$ & 
$y^2\der{x}+x^3y^{-1}\der{y}$ & \\
$C$ &  $\urn{-1}{2}{2}{-1}$ & $x'=y^2,~y'=x^2$ & $y^2\der{x}+x^2\der{y}$ &
Dixonian \\
$D$ & $\urn{-1}{3}{3}{-1}$ & $x'=y^3,~y'=x^3$ & $y^3\der{x}+x^3\der{y}$
& lemniscatic~\cite{FlGaPe05} \\
$E$ & $\urn{-1}{3}{5}{-3}$ & $x'=y^3,~y'=x^5y^{-2}$ & 
$y^3\der{x}+x^5y^{-2}\der{y}$ & \\
$F$ & $\urn{-1}{4}{5}{-2}$ & $x'=y^4,~y'=x^5y^{-1}$ & 
$y^4\der{x}+x^5y^{-1}\der{y}$ & \\
\hline
\end{tabular}
\end{center}

\caption{\label{sixell-fig} The six elliptic cases of a $2\times2$
P\'olya urn: matrix ($\cal M$), system ($\Sigma$), operator
($\delta$), 
and type, following~\cite{FlGaPe05}.}
\end{figure}
\begin{figure}
\begin{center}
\fbox{\Img{11}{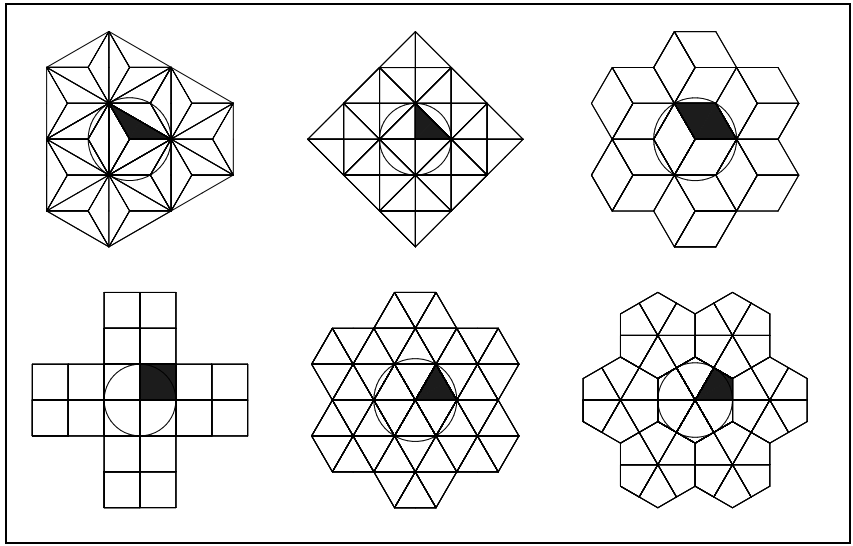}}
\end{center}

\caption{\label{sixlat-fig}
The six lattices $\binom{A~B~C}{D~E~F}$
corresponding to urn models that are exactly
solvable in terms of elliptic functions.}
\end{figure}
The bivariate generating function of urn histories 
(a prototype is provided by Proposition~\ref{alltimes-prop} under the
$\M_{12}$ scenario)
is in each case
expressed in terms of a function~$\psi$ that is determined as the inverse
of Abelian integral over the Fermat curve~$\F_p$ of degree
$p=s+a+b$ (so that the genus is
$(p-1)(p-2)/2$). Reference~\cite{FlGaPe05} starts by building the
bivariate generating function of histories as
$H:=e^{z\delta}[x]$; then, the method of
characteristics is applied in order to solve a partial differential
equation satisfied by~$H$. Elementary conformal mapping arguments eventually
lead to a complete characterization of all urn models
and operators~(\ref{gam}) 
such that $H$ is expressible in terms of elliptic
functions.
	(Caveat: In~\cite{FlGaPe05}, a model is said to be solvable by 
	elliptic functions if its ``base'' function $\psi$ is
	a power, possibly \emph{fractional}, of an elliptic function. 
	Under this terminology, some of the models  eventually turn
	out to be elliptic, 
	though they are a priori associated to Fermat curves of
	genus $>1$.). It is found in~\cite{FlGaPe05} that there are
	altogether only six 
possibilities that correspond to regular tessellations of the Euclidean
plane; see Figures~\ref{sixell-fig} and~\ref{sixlat-fig}.

This classification nicely complements some of Dumont's researches in
the 1980s; 
see~\cite{Dumont86}. In particular Dumont developed a wealth of 
combinatorial connections between the Jacobian elliptic
function, Schett's operator (which involves three variables
rather than two), and permutations. The Dumont--Schett's result should in
particular provide a complete analytic model for the urn with balls of
\emph{three} coulours,
\[
\cal M=\left(\begin{array}{ccc}
-1&1&1\\1&-1&1\\1&1&-1
\end{array}\right),
\]
corresponding to the trivariate operator introduced by Schett~\cite{Schett76}
\begin{equation}\label{schett}
\delta=yz\der{x}+zx\der{y}+xy\der{z}.
\end{equation}
Observe that very few models are known to be explicitly solvable
in the case of  urns with
three types of balls or more\footnote{%
	Puyhaubert's thesis~\cite{Puyhaubert05} provides a
	a discussion of $3\times3$ triangular cases.
	On the other hand, probabilistic techniques 
	yielding asymptotic
	information are available 
	for urns with an arbitrary number of colours: see for instance
	recent studies by Janson and
	Pouyanne~\cite{Janson04a,Pouyanne05}. In particular Pouyanne
	(following 
	a suggestion of Janson) obtains results relative 
	to the Yule process with $k$ types of
	balls. The rules are of the form
	$
	\x_1\lora\x_2\x_2,~ \x_2\lora\x_3\x_3,~\ldots,
	\x_{k-1}\lora\x_{k}\x_{k},~ \x_k\lora \x_1\x_1,
	$
	and a curious ``phase transition'' is found to occur at $k=9$.
}.

The proofs obtained in the present
paper rely on the basic combinatorial properties of 
an associated non-linear differential system. They are
purely ``conceptual'', thereby bypassing several 
computational steps of~\cite{FlGaPe05}, in particular the method of
characteristics. 
The process based on ordinary differential
systems that has been developed here for the $\M_{12}$ urn
is in fact applicable to all urns of type~(\ref{Mgen})
that are balanced ($\alpha+\beta=\gamma+\delta$).
This observation suggests,
more generally,  that interesting combinatorics is likely to be found
amongst several nonlinear autonomous systems. 
Perhaps something along the lines of Leroux and Viennot's
combinatorial-differential calculus~\cite{LeVi88b} is doable here.

\section{Second model: permutations and parity of levels}\label{firstperm-sec}

Our \emph{second combinatorial model} of Dixonian functions is in terms of permutations.
It necessitates the notion of level of an element (a value) in permutation, itself 
related to a basic tree representation, as well as a basic classification 
of elements into four local ordinal types (peaks, valleys, double
rises, and double falls). 

\smallskip

\begin{figure}
\Img{5.5}{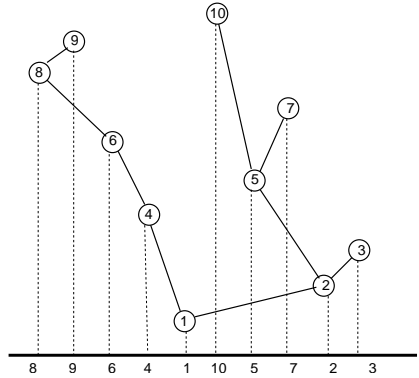}
\caption{\label{bintree-fig}
A permutation of length~10 and its
corresponding increasing binary tree.}
\end{figure}

A permutation can always be represented as a tree, which is binary,
rooted, and increasing
(see Stanley's book~\cite[p.~23]{Stanley86}).
Precisely, let $w=w_1w_2\cdots w_n$ be a word on $\Z_{>0}$ without repeated
letters. A tree denoted by $\Tree(w)$ is associated to a word~$w$ by the following rules.
\begin{itemize}
\item[---] If $w$ is the empty word, then $\Tree(w)=\varepsilon$ is the
empty tree.
\item[---] Else, let $\xi=\min(w)$ be the least element of~$w$. Factor $w$ as
$w' \xi w''$. The tree $\Tree(w)$  is then inductively defined by
\[
\Tree(w)=\langle  \Tree(w'),\xi, \Tree(w'')\rangle,
\]
that is, the root is $\xi$, and $\Tree(w'),\Tree(w'')$ are
respectively the left and right subtrees of the root.
\end{itemize}
The tree so obtained\footnote{%
	 In~\cite[p.~41]{Stanley86}, Stanley
	writes that such models of permutations have been 
	``extensively developed primarily by the French'' and refers
	to Foata-Sch\"utzenberger~\cite{FoSc70}.
} is such that the smallest letter of the word appears
at the root, and the labels of any branch stemming from the root go in
increasing order. This construction applies in particular to any permutation,
$\sigma=\left(\begin{array}{cccc}
1&2&\cdots &n \\
\sigma_1&\sigma_2&\cdots &\sigma_n \end{array}\right)$,
once it is written as the equivalent word $\sigma_1\sigma_2\cdots \sigma_n$.
An infix order traversal of the tree, a projection really, gives back
the permutation from the increasing tree, so that the correspondence is
bijective; see Figure~\ref{bintree-fig}. Observe also that 
we are dealing with  binary trees in the usual sense of computer science~\cite[\S2.3]{Knuth97},
a nonempty node being of one of four types: binary, nullary (leaf), 
left-branching, and right-branching.

\begin{definition}\label{lev-def}
Given a rooted tree, the \emph{level} of a node~$\nu$ is defined as the
distance, measured in the number of edges, between~$\nu$ and the
root. The level of an element (letter) in a word comprised
of distinct letters is
the level of the corresponding node in the associated tree.
\end{definition}

The statement below also makes use of 
a fundamental classification of local order types in permutations;
see in particular Fran{\c c}on and Viennot's paper~\cite{FrVi79}.

\begin{definition}\label{types-def}
Let a permutation be written as a word, $\sigma=\sigma_1\cdots \sigma_n$.
Any value $\sigma_j$ can be classified into four
local ordinal types called
\emph{peaks, valleys, double rises}, and \emph{double falls},
depending on its ordinal relation to its neighbours:
\begin{center}\small
\begin{tabular}{cccc}
\hline
Peaks & Valleys & Double rises & Double falls\\
\hline
$\sigma_{j-1}<\sigma_j>\sigma_{j+1}$ &
$\sigma_{j-1}>\sigma_j<\sigma_{j+1}$ &
$\sigma_{j-1}<\sigma_j<\sigma_{j+1}$ &
$\sigma_{j-1}>\sigma_j>\sigma_{j+1}$ \\
\hline
\end{tabular}
\end{center}
\end{definition}
For definiteness, \emph{border conditions} 
conditions must be adopted,
and we shall normally opt for one of the two choices,
\[
\sigma_0=-\infty,\qquad
\sigma_{n+1}\in\{-\infty,+\infty\}.
\]
(For instance, an alternating permutation 
is characterized by the fact that it has only
peaks and valleys.)
The
corresponding quadrivariate exponential generating function 
of all permutations was first determined by Carlitz.

\begin{theorem}[Second combinatorial model]\label{parity-prop}
Consider the class $\cal X$ of permutations bordered by
$(-\infty,-\infty)$
such that
elements at any odd level are valleys  only. Then the exponential
generating function is
\[
X(z)=\smh(z)=-\sm(-z).\]
For the class $\cal Y$ of permutations also bordered by $(-\infty,-\infty)$
such that
elements at any even level are valleys  only, the exponential
generating function is
\[
Y(z)=\cmh(z)=-\cm(-z).\]
\end{theorem}
\begin{proof} 
This statement\footnote{%
	As we learnt in May 2005, this result was also obtained 
	independently by
	Dumont at Ouagadougou in 1988; see 
	his unpublished note~\cite{Dumont88}.}
 easily results from 
first principles of combinatorial analysis.
The interpretation is obtained by examining
the differential system
\[
X'(z)=Y(z)^2, \quad Y'(z)=X(z)^2,\qquad
X(0)=0, \quad Y(0)=1,
\]
or, under an equivalent integral form,
\[
X(z)=\int_0^z Y(w)^2\, dw, \qquad
Y(z)=1+\int_0^z X(w)^2\, dw.\]
It suffices to note that, if $\cal A$ and $\cal B$ are two
combinatorial classes with exponential generating functions $A$ and $B$,
then the product $A\cdot B$ enumerates the labelled product $\cal
A\star \cal B$. Also, the integral $\int A\cdot B$  enumerates
all well-labelled triples $\langle \alpha,\xi, \beta\rangle$, where
$\xi$ is the smallest of all labels. 

By the preceding remarks, we see that $X,Y$ are the generating
function of increasing trees satisfying the relations
\[
\cal X=\langle \cal Y,\min, \cal Y\rangle, \qquad
\cal Y = \varepsilon + \langle \cal X,\min, \cal X\rangle.
\]
It is then apparent that trees of type~$\cal X$ are such that they
only have double nodes at odd levels,
while those of type~$\cal Y$ only have double nodes at even levels.
The proof of the proposition concludes given the 
correspondence
between tree node and permutation value types, to wit,
\begin{center}\def\la{$\leftrightarrow$}\renewcommand{\tabcolsep}{3truept}
\begin{tabular}{rclcrcl}
double node &\la& valley , &&
left-branching node &\la& double fall,  \\
leaf &\la&  peak,&&
right-branching node &\la& double rise,
\end{tabular}
\end{center}
which is classical (and obvious via projection).
\end{proof}

Consequently, for $n=3\nu$, the number of permutations of type $\cal
Y$ is $Y_n=n![z^n]\cmh(z)$. We have
\[
\begin{array}{lll}
Y_0=1 && \cal Y_2=\{\epsilon\}\quad\hbox{(the empty permutation)}\\
Y_3=2 && \cal Y_2 = \{213,~312\},
\end{array}
\]
which agree with~(\ref{sctayl}). In order to form the shapes of trees
of type~$\cal Y$, one can use the grammar represented graphically as
($\Box$ represents the empty tree)
\[
\cal Y\quad = \quad \Box \quad + \quad  
\hbox{\setlength{\unitlength}{1.truemm}\thicklines
\begin{picture}(24,16)
\put(12,0){\circle*{2}}
\put(6,6){\circle*{2}}
\put(6,6){\line(1,-1){6}}
\put(18,6){\circle*{2}}
\put(18,6){\line(-1,-1){6}}
\put(1,14){$\cal Y$}
\put(2,12){\line(2,-3){4}}
\put(9,14){$\cal Y$}
\put(10,12){\line(-2,-3){4}}
\put(13,14){$\cal Y$}
\put(14,12){\line(2,-3){4}}
\put(21,14){$\cal Y$}
\put(22,12){\line(-2,-3){4}}
\end{picture}}\]
\begin{figure}
\begin{center}
\def\b{\hspace*{0.35truecm}}
\Img{6}{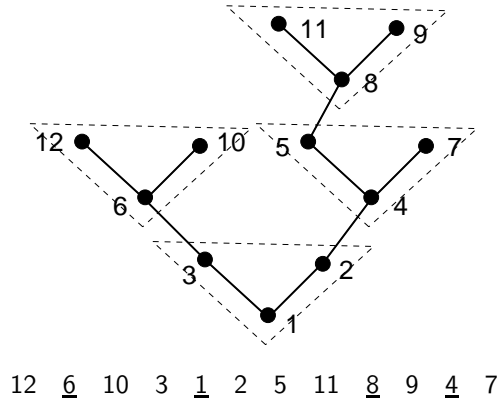}
\\[2truemm]
{\sf 12\b \underline{6}\b 10\b 3\b \underline{1}\b 2\b 5\b 11\b
\underline{8}\b 9\b \underline{4}\b 7}
\end{center}
\caption{\label{cmperm-fig}
The construction of a $\cal Y$-permutation enumerated by $\cm(z)$. The elements
at even level (underlined) are valleys only.}
\end{figure}
This automatically generates trees whose sizes (their number of nonempty nodes) are
multiples of~3. (The number of such tree shapes of size $n=3\nu$ is the same as
the number of quaternary trees of size~$\nu$, namely,
$\frac{1}{3\nu+1}\binom{4\nu}{\nu}$). The increasing labellings of
these trees in all possible ways followed by projections provides a
way of listing all permutations of type~$\cal Y$. (See Figure~\ref{cmperm-fig}
for an illustration of this construction.) In this way,
it is found that there are~4
 legal tree shapes of type~$\cal Y$, each of which admitting of
10 increasing labellings, which globally corresponds to $Y_6=40$, 
in agreement with~(\ref{sctayl}). 

Similarly, trees corresponding to $\cal X$ are determined by the
grammar
\[
\cal X\quad = \quad 
\hbox{\setlength{\unitlength}{1.truemm}\thicklines
\begin{picture}(24,9)
\put(12,0){\circle*{2}}
\put(4,7){$\cal Y$}
\put(6,6){\line(1,-1){6}}
\put(17,7){$\cal Y$}
\put(18,6){\line(-1,-1){6}}
\end{picture}}\]
which yields, again in agreement with~(\ref{sctayl}):
\[\begin{array}{lll}
X_1=1 : && \cal X_1=\{1\}, \\
X_4=4 : && \cal X_4=\{3241,~4231,~1324,~1423\}.
\end{array}
\]

\begin{Note} \emph{Parity-based permutation models for elliptic
functions.} \label{parperm-note}
There already exist several combinatorial models of the Jacobian elliptic
functions $sn,cn$ due to Viennot, Flajolet and Dumont\footnote{%
	These authors were answering a question of Marco
	Sc\"hutzenberger in the 1970's.
	Sch\"utzenberger first conjectured that Jacobian elliptic
	functions should have combinatorial content since their coefficients
	involve both factorial and Euler numbers.
} 
and discovered around 1980.
These all
involve permutations restricted by a parity condition of sorts.
\begin{itemize}
\item[---] Viennot~\cite{Viennot80} develops a model based on the
differential system satisfied by the Jacobian elliptic functions,
namely
\[
x'=yz, \qquad y'=-zx,\qquad z'=-k^2xy,
\]
where $x,y,z$ represent the classical Jacobian functions $sn,cn,dn$. 
This leads to an interpretation of the coefficients in terms of 
a class of permutations called by Viennot ``Jacobi
permutations''. Such permutations satisfy \emph{parity restrictions}
(mimicking the differential system) and are enumerated by Euler numbers.
\item[---] Flajolet~\cite{Flajolet80b} observed, from continued
fraction theory, that the quantity
\[
(-1)^n(2n)!\cdot [z^{2n}\alpha^{2k}]\, cn(z,\alpha)
\]
counts alternating
permutations of length $2n$ that have $k$ valleys of \emph{even value}.
\item[---] Dumont~\cite{Dumont79,Dumont81} provides an elegant
interpretation of the coefficients of $sn,cn,dn$ in terms of the
\emph{parity} of peaks of cycles in the cycle decomposition of
permutations. Dumont's results are based on consideration of Schett's
partial differential operator, already encountered in~(\ref{schett}).
\end{itemize}
Theorem~\ref{parity-prop} adds another parity-based permutation model
to the list.
\end{Note}

\section{Third model: permutations and repeated patterns}\label{secondperm-sec}

Our \emph{third  combinatorial model} of
Dixonian functions is again
in terms of permutations. It 
relies on repeated permutations, much
in the style of Flajolet and Fran{\c c}on's earlier interpretation of Jacobian
($sn,cn$) elliptic functions~\cite{FlFr89}, but different. It is
notable that the 
Flajolet-Fran{\c c}on permutations are based on a binary pattern, whereas
those needed here involve a symmetry of order~three.
The type of an element (a value) in a permutation is, as in the
previous section, any of the set peak, valley, double rise,
double fall.
\begin{definition} An \emph{$r$--repeated} permutation 
of size~$n$ is a permutation
such that for each~$j$ with $j\ge 0$, the elements of \emph{value} in
$\{jr+1,jr+2,\ldots,jr+r-1\}\cap\{1,\ldots,n\}$ are all of the same
ordinal type, namely all peaks, valleys, double rises, or double falls.
\end{definition}
For instance, the permutation
\def\pp{\hbox{\tiny${\nearrow\!\searrow}$}}
\def\vv{\hbox{\tiny${\searrow\!\nearrow}$}}
\def\rr{\hbox{\tiny${\nearrow\!\!\nearrow}$}}
\def\ff{\hbox{\tiny${\searrow\!\!\searrow}$}}
\begin{center}\sf
\renewcommand{\tabcolsep}{2.5truept}
\begin{tabular}{ccccccccccccc}
12&8&7&2&5&10&3&6&13&9&1&4&11\\
\pp&\ff&\ff&\vv&\rr&\pp&\vv&\rr&\pp&\ff&\vv&\rr&\pp
\end{tabular}
\end{center}
bordered with $(-\infty,-\infty)$ is a 3-repeated permutation of
size~13, as is immediately verified by the listing of value types
on the second line.

Dixonian functions will be proved to enumerate a variety of 3--repeated
permutations. The proof is indirect and it first necessitates
Flajolet's combinatorial theory of continued
fractions~\cite{Flajolet80b,GoJa83} as well as a bijection between
a system of weighted lattice paths and permutations, of which a first
instance was discovered by Fran{\c c}on and Viennot in~\cite{FrVi79}. 

\subsection{Combinatorial aspects of continued fractions.}\label{confrac-sec}
Define a \emph{lattice path},
also known as a \emph{Motzkin path}, of length~$n$ as a sequence 
of numbers $s=(s_0,s_1,\ldots,s_n)$, satisfying the conditions
\[
s_0=s_n=0,\quad s_j\in \Z_{\ge0}, \qquad
\left|s_{j+1}-s_j\right|\in\{-1,0,+1\}.\]
This is represented as a polygonal line in the Cartesian plane
$\Z\times\Z$. A step is an edge $(s_j,s_{j+1})$, and it is said to be
an ascent,  
a level, or a descent  according to the value, respectively~$+1,0,-1$, of
$s_{j+1}-s_j$; the quantity $s_j$ is called the (starting) altitude of
the step. A path without level steps is a \emph{Dyck path}.
Motzkin paths are enumerated by Motzkin numbers; Dyck paths belong to
the Catalan realm~\cite{Stanley86}.

Let $P({\bf a},{\bf b},{\bf c})$ be the  generating
function of lattice paths in 
infinitely many indeterminates ${\bf a}=(a_k)$, ${\bf b}=(b_k)$,
${\bf c}=(c_k)$, with $a_k$ marking an ascent from
altitude~$k$ and similarly for descents marked by $b_k$ and for levels
marked by~$c_k$. In other words, associate to each path $\varpi$
a monomial $\frak{m}(\varpi)$ like in this example:
\[\renewcommand{\tabcolsep}{3truept}
\hbox{\begin{tabular}{ccc}$\varpi={}$ & \begin{tabular}{c}
\def\circo{\circle*{2}}\setlength{\unitlength}{0.5truemm}\thicklines
\begin{picture}(100,30)
\put(0,0){\circo}
\put(0,0){\line(1,1){10}}\put(10,10){\circo}\put(10,10){\line(1,1){10}}
\put(20,20){\circo}\put(20,20){\line(1,0){10}}\put(30,20){\circo}
\put(30,20){\line(1,-1){10}}\put(40,10){\circo}\put(40,10){\line(1,1){10}}
\put(50,20){\circo}\put(50,20){\line(1,-1){10}}\put(60,10){\circo}
\put(60,10){\line(1,-1){10}}\put(70,0){\circo}\put(70,0){\line(1,1){10}}
\put(80,10){\circo}\put(80,10){\line(1,-1){10}}\put(90,0){\circo}
\thinlines\put(0,0){\line(1,0){90}}
\end{picture}
\end{tabular}&\begin{tabular}{c} $\begin{array}{lll}\frak{m}(\varpi)
&=&a_0a_1c_2b_2a_1a_2b_3b_2b_1a_0b_1 \\
&=&a_0^2a_1^2a_2b_1^2b_2^2b_3c_1;
\end{array}$
\end{tabular}\end{tabular}}
\]
then, one has
\[P({\bf a},{\bf b},{\bf c})=\sum_{\varpi} \frak{m}(\varpi),\]
where the sum is extended over all Motzkin paths~$\varpi$.

Here is what Foata once called ``the shallow Flajolet
Theorem''\footnote{%
	This designation stands to reason
	as the proof is extremely easy, so that
	the theorem borders on being an ``observation''. However,
	the paper~\cite{Flajolet80b} is really a 
	``framework'' where orthogonal polynomial systems, lattice
	paths, continued fractions, Hankel determinants, etc, all find
	a combinatorial niche.
}  taken from~\cite{Flajolet80b}: 
 
\begin{theorem}[Flajolet~\cite{Flajolet80b}] \label{cf-thm}
The generating function in infinitely many variables
enumerating all lattice paths according to ascents, descents, levels,
and corresponding altitudes is
\begin{equation}\label{flajcf}
P({\bf a},{\bf b},{\bf
c})=\cfrac{1}{1-c_0-\cfrac{a_0b_1}{1-c_1-\cfrac{a_1b_2}{1-c_2-\cfrac{a_2b_3}{\ddots}}}}.
\end{equation}
\end{theorem}
\begin{proof} All there is to it is the correspondence between 
a rational generating function given by a quasi-inverse and certain
Motzkin paths of height $\leq0$. Figuratively:
\[\def\circo{\circle*{2}}\setlength{\unitlength}{0.6truemm}\thicklines
\begin{array}{ccc}
\ds \frac{1}{\ds 1-c_0-a_0Y_1b_1} 
& \equiv & 
\left\{\begin{array}{c}\hbox{~~\begin{picture}(56.5,15)(-25,0)
\put(-25,0){\circo}
\put(-25,0){\line(1,0){10}}
\put(-22,2){$c_0$}
\put(-15,0){\circo}
\put(-9,0){+}\put(0,0){\circo}\put(0,0){\line(1,1){10}}\put(10,10){\circo}
\put(0.5,7){$a_0$}
\put(10,10){\line(1,0){10}}\put(20,10){\circo}
\put(13,12){$Y_1$}
\put(20,10){\line(1,-1){10}}
\put(24,7){$b_1$}
\put(30,0){\circo}
\end{picture}} 
\end{array}\right\}^{\ds\star}\,. \end{array}
\]
It then suffices to apply repeatedly substitutions of the form
\[
Y_1\mapsto (1-c_1-a_1Y_2b_2)^{-1}, \quad
Y_2\mapsto (1-c_2-a_2Y_3b_3)^{-1}, \cdots\,,
\]
in order to generate the continued fraction of the statement.
\end{proof}

A \emph{bona fide} generating function obtains by making the
generating function homogeneous, setting
\[
a_j\mapsto \alpha_j z, \qquad 
b_j\mapsto \beta_j z, \qquad 
c_j\mapsto \gamma_j z
\]
where $\alpha_j,\beta_j,\gamma_j$ are new formal variables. 
Once these variables are assigned numerical values (for instance
$\alpha_j=\beta_j=\gamma_j=j+1$), the generating function
of~(\ref{flajcf}) enumerates weighted lattice paths by length (marked
by~$z$), with weights $\{\alpha_j\}\cup\{\beta_j\}\cup\{\gamma_j\}$
taken multiplicatively. For ease of reference, we encapsulate this notion into
a formal definition. 

\begin{definition} A system of path diagrams is the class of
multiplicatively weighted Motzkin paths determined by a
possibility function~$\Pi$, 
which assigns numerical values to the formal variables
$\alpha_j,\beta_j,\gamma_j$. 
\end{definition}

The values are normally integers, in which case a particular path
diagram is equivalent to a lattice path augmented by a sequence
of ``choices'' of the same length,
the number of possible choices being $\alpha_j$ for an ascent from
altitude~$j$, and so on. Various
systems of path diagrams are known to be bijectively
associated to
permutations, involutions, set partitions, and preferential arrangements,
to name a few~\cite{Flajolet80b,Flajolet82,FrVi79,GoJa83}.

\begin{Note} \emph{Andr\'e's method
and an alternative derivation of the $J$-fractions
relative to $\sm,\cm$.} \label{andre-note}
D\'esir\'e Andr\'e published in 1877 a remarkable study~\cite{Andre77},
in which he was 
able to represent the coefficients of the Jacobian functions $sn,cn$
as arising from weighted lattice paths, that is, from a system of path
diagrams. We apply here Andr\'e's methodology to Dixonian functions.
A direct use of the continued fraction Theorem~\ref{cf-thm},
then  provides a 
direct derivation of the $J$-fractions relative to $\sm,\cm$,
which relies in a simple way on basic algebraic properties
of the fundamental differential system~(\ref{dixon0}).

Let us consider the case of $\smh$, which is the $s$-component of the
usual system $s'=c^2$, $c'=s^2$. From this system and
elementary algebra induced by $c^3-s^3=1$, we find that $s$ satisfies
a third-order non linear differential equation,
$s'''=6s^4+4s $, and more generally ($\partial$ represents 
differentiation with respect to the independent variable):
\begin{equation}\label{smm}
\partial^3 s^m =\left( m(m+1)(m+2)s^3 +
2m(m^2+1)+\frac{m(m-1)(m-2)}{s^3} \right)\, s^m .
\end{equation}
This shows that there exists a family of polynomials $(P_k)$ with
$\deg(P_k)=3k+1$ such that
\[ \partial^{3k} s = P_k(s),\]
where $P_k$ is of the form $P_k(w)=w\widehat P_k(w^3)$, with $\widehat
P_k$ itself a polynomial.
Then, Taylor's formula provides
\begin{equation}\label{smhtayl}
\smh(z)=\sum_{k=1}^\infty P_k'(0) \frac{z^{3k+1}}{(3k+1)!} 
\end{equation}
upon taking into account the fact that only coefficients of index
$1,4,7,\ldots$ in $\smh$ are nonzero. Thus, the Taylor coefficients of $\smh$
are accessible from the coefficients of the lowest degree monomials in
the $P_k$ polynomials.

Introduce now a notation for the coefficients of the $P_k$ polynomials:
\[
P_{k,m}=[w^m] P_k(w).\]
The basic equation~(\ref{smm}) implies the recurrence ($m\equiv1 \pmod
4$).
\[
P_{k+1,m}=(m+1)(m+2)(m+3)P_{k,m+3}+2m(m^2+1)P_{k,m}+(m-1)(m-2)(m-3)P_{k,m-3}.\]
This relation expresses precisely the fact that the coefficient
$[w^m]P_k(w)$ is the number of weighted paths starting at the point
$(0,1)$ in the lattice $3\Z\oplus (1+3\Z)$ where  the
elementary steps are of the form $\vec a=({1},{+3})$,
$\vec b=({1},{-3})$, $\vec c=(1,{0})$, and the weights of 
steps starting at an ordinate~$m\equiv1\pmod4$ are  respectively
\begin{equation}\label{weights}
\widehat\alpha_m=m(m+1)(m+2),\quad
\widehat\beta_m=m(m-1)(m-2),\quad
\widehat\gamma_m=2m(m^2+1).
\end{equation}
Up to a vertical translation of~$-1$ and a vertical rescaling of the
lattice by a factor of $\frac13$, the coefficient $P_{1,\nu}$
is seen to
enumerate standard weighted lattice paths (having steps $0,\pm1$)
of length~$3\nu$ 
with the new weights for steps starting at altitude~$\ell$ being:
\[
\alpha_\ell=(3\ell+1)(3\ell+2)(3\ell+3),\quad
\beta_\ell=(3\ell-1)(3\ell)(3\ell+1),\quad
\gamma_\ell=2(3\ell+1)((3\ell+1)^2+1).
\]
Theorem~\ref{cf-thm} then yields the continued fraction expansion of
the ordinary generating function associated to~$\smh$
(equivalently, to $\sm$ under a simple change of signs).
Similar calculations apply to the other $J$-fractions of
Section~\ref{conrad-sec}, including the one relative to $\cm$.

It is especially interesting to observe the way in which the original
differential system satisfied by $\smh,\cmh$ churns out weighted lattice
paths by means of higher order differential relations. 
This is the spirit of Andr\'e's work who played the original
game on the Jacobi normal form of elliptic functions taken under the
form (Andr\'e's notations):
\[
\left[\frac{d\varphi(x)}{dx}\right]^2=\cal D + \cal V
\varphi^2(x)+\cal G\varphi^4.\]
Obviously, Andr\'e intends $\cal D$ to mean
``right'', $\cal G$ to mean ``left'' and $\cal V$ to mean
``vertical'', as he had in mind paths rotated by $90^{\circ}$. 
(In one more step, 
Andr\'e could have
discovered the continued fractions associated to $sn,cn$ had he known
Theorem~\ref{cf-thm} or the earlier technique of ``Stieltjes
matrices''~\cite{Stieltjes89},
 whose determination
is equivalent to a continued fraction expansion.)
\end{Note}

\begin{Note} \emph{On addition theorems.} \label{add-note}
A standard way to derive explicit continued fraction
expansions is by means of Rogers' addition theorem (itself logically
equivalent to a diagonalization technique of Stieltjes).
We say that an analytic
function defined near~0
(or a formal power series) $f(z)$ satisfies an \emph{addition
formula} if it can be decomposed in terms of a sequence $\phi_k$ of
functions as
\begin{equation}\label{sradd}
f(z)=\sum_{k\ge0} \phi_k(x)\cdot \phi_k(y),\qquad\hbox{where}\quad
\phi_k(x)\mathop{=}_{x\to 0} O(x^k).
\end{equation}
Then knowledge of the coefficients $[x^k]\phi_k(x)$ and
$[x^{k+1}]\phi_k(x)$ implies knowledge of the coefficients in the
$J$-fraction representation of the formal Laplace transform $F$ of $f$,
$F(s)=\Lap(f,s)$: the  formul{\ae} are simple, see
Theorem~53.1 in Wall's treatise~\cite[p.~203]{Wall48}. 
A good illustration is  
the addition formula for $\sec(z)$, namely
\[
\sec(x+y)=\frac{1}{\cos x\cos y -\sin x\sin y}=
\sum_{k\ge0} \sec x \tan^k x \cdot \sec y  \tan^k y,
\]
which provides analytically the continued 
fraction expansion of the ordinary
generating function of Euler numbers derived
combinatorially below as~(\ref{tansec-cf}).

Dixonian functions,
being elliptic functions, are known to admit addition formul{\ae},
albeit of a form different from~(\ref{sradd}). For instance, one has
\[
\begin{array}{lll}
\cm(u+v)& = &\ds \frac{c_1c_2-s_1s_2(s_1c_2^2+s_2c_1^2)}{1-s_1^3s_2^3}\\
&& \hbox{with}~s_1=\sm(u),~c_1=\cm(u),~s_2=\sm(v),~c_2=\cm(v),
\end{array}
\]
found in Dixon's paper (\S33, Equation~(39), p.~183
of~\cite{Dixon90}). It is in fact possible to combine this elliptic addition
theorem with developments from the previous note in order to 
come up with an addition theorem of Stieltjes-Rogers
type~(\ref{sradd}) relative to Dixonian
functions: the form of the $\phi_k$ is essentially given by
$P_k(\sm(z))$  and the generating function for the $P_k$ polynomials
can be explicitly determined (details omitted).
\end{Note}

\subsection{Correspondences between lattices paths and permutations.}
We next need bijections due to
 Fran{\c c}on-Viennot~\cite{FrVi79} between
permutations and two systems of path diagrams.

It is convenient to start the discussion by introducing what
V.I.~Arnold~\cite{Arnold00} calls \emph{snakes}. Consider piecewise monotonic
smooth functions from $\R$ to $\R$, such that all their critical
values (i.e., the values of their maxima and minima)
are different, and take the
equivalence classes up to orientation preserving
maps of $\R\times\R$. It is sufficient to restrict attention to the two types,
\begin{equation}\label{snaketypes}
\begin{array}{c}\hbox{\Img{7}{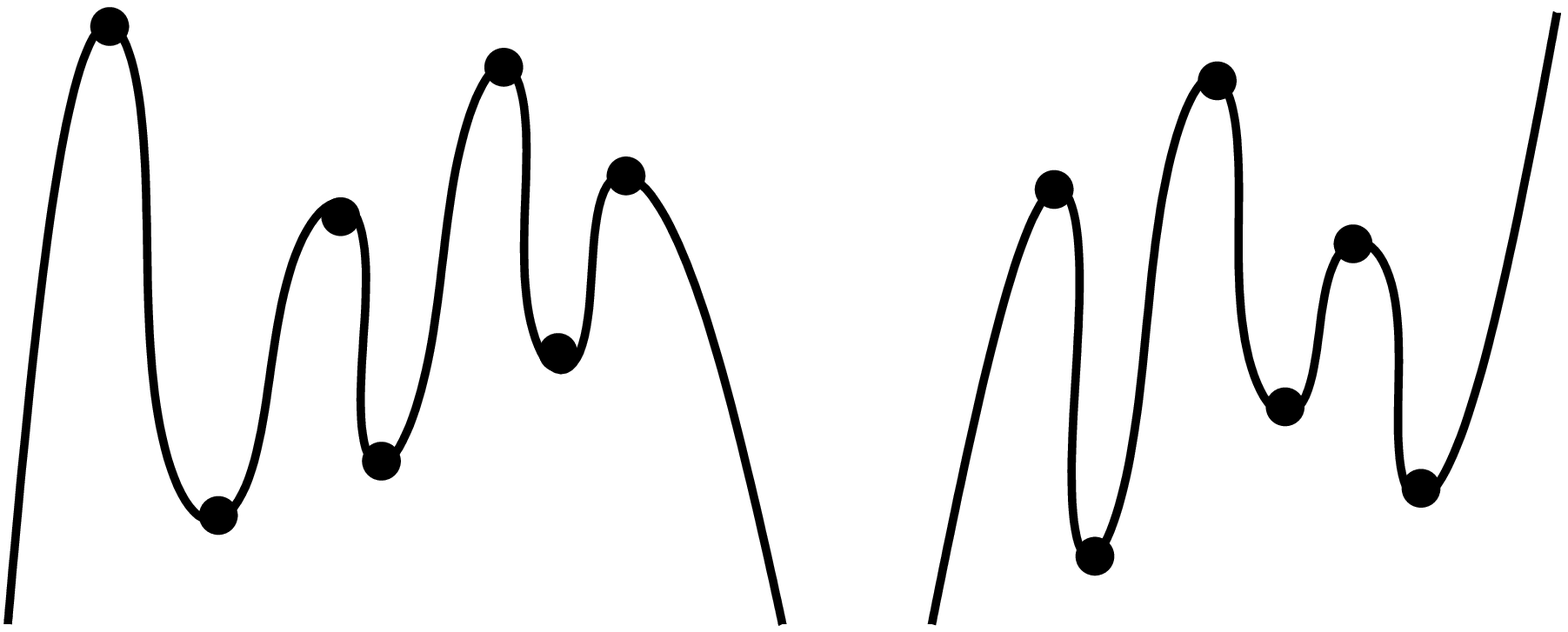}}\\
\begin{array}{ccccccc}(7&1&4&2&6&3&5)\end{array}
\qquad
\begin{array}{cccccc}(5&1&6&3&4&2)\end{array}
\end{array}
\end{equation}
respectively called the $(-\infty,-\infty)$ and $(-\infty,+\infty)$
types. Clearly an equivalence class is an alternating permutation
(see the second line of~(\ref{snaketypes})). The
exponential generating functions corresponding to the two types
of~(\ref{snaketypes}) are by Andr\'e's classic theorem respectively,
\[
\tan(z)=\frac{\sin(z)}{\cos(z)},\qquad \sec(z)=\frac{1}{\cos(z)},
\]
size being measured by the number of critical points and equivalently
the length of the corresponding permutation. 

\begin{figure}
\begin{center}
\Img{9}{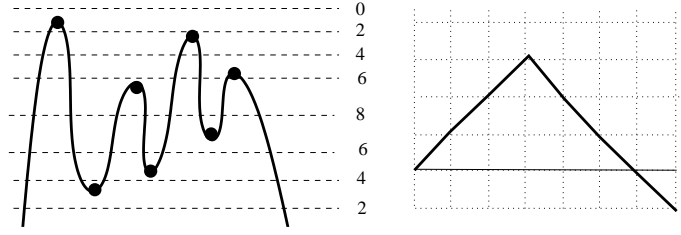}
\end{center}
\caption{\label{sweep-fig}
The sweepline algorithm: a snake and its associated Dyck path.}
\end{figure}

A simple sweepline algorithm associates to a snake a Dyck path
as follows. Consider first the $(-\infty,-\infty)$ case.
Imagine  moving a line from $-\infty$ in the vertical
direction towards $+\infty$. At ordinates that are a midpoint between
two successive critical values, associate a number which is the number
of intersection points of the line with the curve. For a snake with
$2n+1$ critical points, this gives us a
sequence of numbers $x=(x_0,x_1,\ldots,x_{2n+2})$ such that
$x_{j+1}-x_j=\pm2$, $x_0=2$, and $x_{2n+1}=0$. The rescaled sequence
\[
\xi=(\xi_0,\xi_1,\ldots,\xi_{2n+1}), \qquad
\xi_j:=(x_j-2)/2,\]
becomes a standard Dyck path (Figure~\ref{sweep-fig}), which obviously
depends only 
on the underlying alternating permutation.

Alternating permutations are far more numerous than Dyck
paths, so that they must be supplemented by additional
information in order to obtain a proper encoding.
To this effect, introduce 
the system of path diagrams given by the possibility
rule
\begin{equation}\label{oddalt}
\Pi^{\operatorname{odd}} :
\Pi^{(-\infty,-\infty)} :
\qquad
\alpha_j =(j+1),\quad \beta_j=(j+1),\quad \gamma_j=0,
\end{equation}
which, we claim, is now bijectively associated to odd-length alternating
permutations,
\[
-\infty<\sigma_1>\sigma_2<\cdots <\sigma_{2m+1}>-\infty\,.
\]
This is easily understood as follows. When executing the
sweepline algorithm, there are several places (possibilities)
 at which a local
maximum (or peak) and a local minimum (or valley) 
can be attached by ``capping'' or ``cupping'':
\[
\hbox{\Img{11}{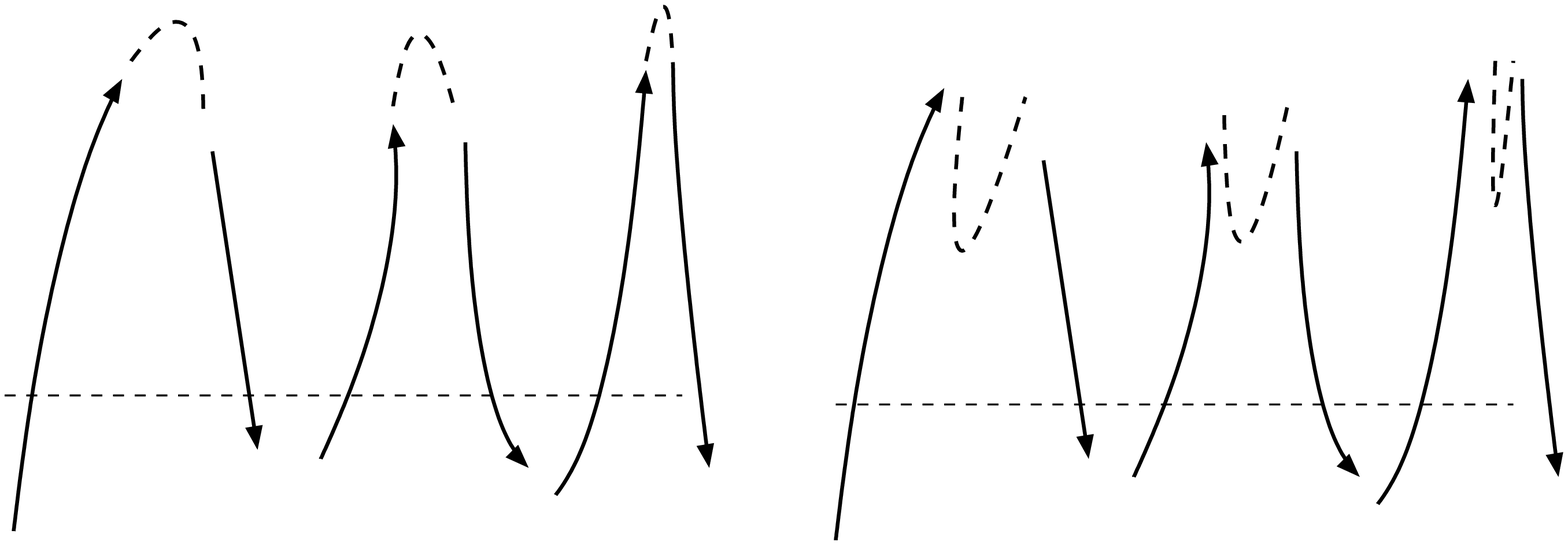}}
\]
These possibilities, as a function of altitude (the
$\xi_j$) are seen  to correspond 
exactly to the possibility set~(\ref{oddalt}).
(For instance, on the example, the rescaled Dyck path is at altitude 2,
there are three possibilities for cupping and three possibilities for
capping.)

For the $(-\infty,+\infty)$ case, that is, even-length alternating
permutations, a similar reasoning shows that the possibilities are
\begin{equation}\label{evenalt}
\Pi^{(-\infty,+\infty)} :
\qquad
\alpha_j =(j+1),\quad \beta_j=j,\quad \gamma_j=0,
\end{equation}
See:
\[
\hbox{\Img{11}{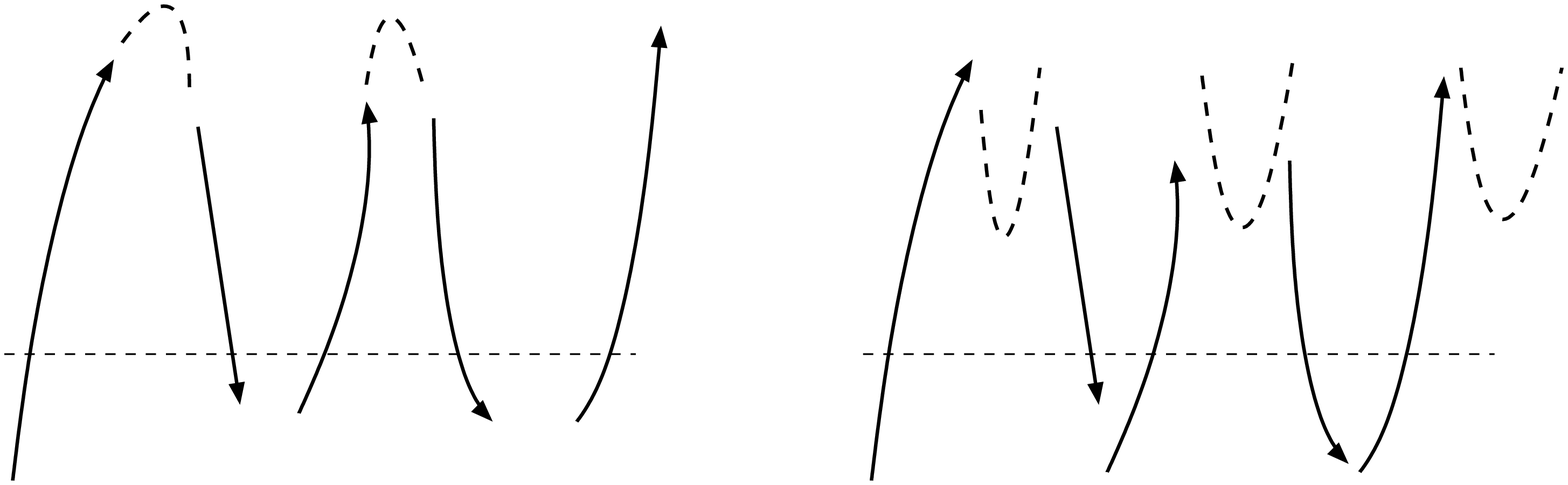}}
\]

By the continued fraction theorem, these two encoding yield 
two continued
fraction expansions originally discovered by Stieltjes (the Laplace
transforms are taken as ``formal'' here):
\begin{equation}\label{tansec-cf}
\left\{\begin{array}{lll}
\ds \int_0^\infty 
\tan(zt) e^{-t}\, dt &=& \ds
\cfrac{z}{1-\cfrac{1\cdot 2\, z^2}{1-\cfrac{2\cdot 3\,
z^2}{\ddots}}},
\\
\ds \int_0^\infty 
\sec(zt) e^{-t}\, dt &=& \ds
\cfrac{1}{1-\cfrac{1^2\, z^2}{1-\cfrac{2^2\,
z^2}{\ddots}}}.
\end{array}\right.
\end{equation}

The previous two bijections can be modified so as to take
into account \emph{all} permutations, not just alternating ones.
What this corresponds to is \emph{spotted snakes}, which are snakes
augmented with an arbitrary finite number of 
non-critical points that are distinguished. (As usual, one operates up to
topological equivalence and the spotted points must have different
altitudes.) 
It then suffices to encode nodes on upward and downward slopes
by level steps of the system of path diagrams in order to get Motzkin
paths~\cite{Flajolet80b,FrVi79,GoJa83}.
The path diagrams so obtained have closely resembling possibility rules:
\[
\left\{\begin{array}{llll}
\overline{\Pi}^{(-\infty,-\infty)} :
&
\alpha_j =(j+1),& \beta_j=(j+1),& \gamma_j=2j+2
\\
\overline{\Pi}^{(-\infty,+\infty)} :
&
\alpha_j =(j+1),& \beta_j=j,& \gamma_j=2j+1.
\end{array}\right.
\]
(In the $(-\infty,+\infty)$ case, one can never insert a cap or a descending node on
the extreme right,
 so that one possibility is suppressed for descents and level steps
of the path diagrams.) The resulting continued fractions are then
\begin{equation}\label{fac-cf}
\left\{\begin{array}{lll}
\ds \sum_{n=0}^\infty (n+1)! z^{n+1}
&=&\ds \cfrac{z}{1-2z-\cfrac{1\cdot 2\, z^2}{1-4z-\cfrac{2\cdot 3\,
z^2}{\ddots}}},
\\
\ds \sum_{n=0}^\infty n!z^n
&=&\ds \cfrac{1}{1-z-\cfrac{1^2\, z^2}{1-3z-\cfrac{2^2\,
z^2}{\ddots}}}.
\end{array}\right.
\end{equation}
These last two fractions are originally due to Euler. 
(The last four continued fraction
expansions were first established  combinatorially by Flajolet
in~\cite{Flajolet80b}.) 

\begin{figure}

\begin{center}
\Img{6}{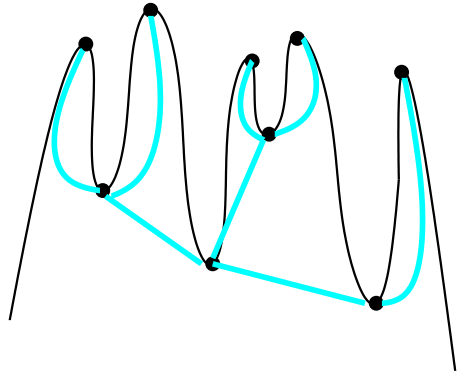}
\end{center}

\caption{\label{sweepviennot-fig} A figure suggesting that the
Fran{\c c}on--Viennot correspondence and the sweepline algorithms are
one and the same thing.}
\end{figure}

\begin{Note} \emph{Snakes of bounded width.}
The usual enumeration of snakes up to deformation is usually presented
by the Russian School as resulting from the ``Seidel-Entringer-Arnold
triangle''. The continued fraction connection exposed here gives access to
new parameters, and, in particular, the ones associated with
convergents of a basic continued fraction.
In this context, it provides the
ordinary generating function of odd snakes of bounded \emph{width},
where width is defined as the maximal cardinality of the image of any
value:
\def\width{\operatorname{width}}\def\card{\operatorname{card}}
\[
\width(s):=\max_{y\in \R} \card\left\{ x~\bigm|~s(x)=y\right\}.
\]
(Width is a trivial variant
of the clustering index introduced in~\cite[p.~159]{Flajolet80b}.)
Take for definiteness snakes with boundary condition $(-\infty,+\infty)$.
Transfer matrix methods first imply that 
the generating function $W^{[h]}(z)$ of snakes having width at most~$2h-1$ 
is \emph{a priori} a rational function,
\[
W^{[h]}(z)=\frac{\overline{P}_h(z)}{\overline{Q}_h(z)},
\]
the first few values being, for width at most 1, 3, 5, 7,
\[
\frac{1}{1},\quad \frac{1}{1-z^2},\quad
\frac{1-4z^2}{1-5z^2},\quad\frac{1-13z^2}{1-14z^2+9z^4}.
\]
(For instance, width~$\leq3$ corresponds to the type $(-\infty,2,1,4,3,6,5\ldots,+\infty)$.)
Continued fraction theory~\cite{Flajolet80b}
further implies that $W^{[h]}(z)$ is
a convergent of the continued fraction relative to $\sec(z)$
in~(\ref{tansec-cf}), with the denominators $\overline Q_h(z)$
being reciprocal polynomials  of
an orthogonal polynomial system. Here, one finds
\[
\overline{Q}_h(z)=z^h Q_h(1/z),\qquad\hbox{where}\quad
Q_h(z)=[t^h](1+t^2)^{-1/2}\exp(z\arctan t),\]
and the $Q_h$ are the Meixner polynomials~\cite{Chihara78}.
(This can be verified via the generating function 
of the $Q_h$, which satisfies a
differential equation of order~1; the method also  adapts to the
determination of the associated $P_h$ polynomials.)
\end{Note}

\begin{Note} \emph{On Fran{\c c}on-Viennot.} \label{frvi-note}
The bijective encoding of permutations by path diagrams as presented
here is exactly the
same as the one obtained from the original Fran{\c
c}on-Viennot
correspondence---only our more geometric presentation differs. 
Indeed, one may think of the Fran{\c c}on--Viennot
correspondence given in~\cite{FrVi79} as the gradual construction of an increasing
binary tree, upon successively appending
nodes at dangling links~\cite{Flajolet80b,GoJa83}; see also
Figure~\ref{sweepviennot-fig},
where both a snake and the underlying tree are represented.
Use will made below of the tree view of the Fran{\c c}on-Viennot correspondence.
(Biane~\cite{Biane93}
discovered a related sweepline algorithm, but 
it is applied to the decomposition
of permutations into cycles.)
\end{Note}

\subsection{The models of r--repeated permutations.}
The Fran{\c c}on-Viennot bijection and the Continued Fraction
Theorem provide:

\begin{proposition}\label{r-prop}
Let $R_{rn+1}$ be the number of $r$--repeated permutations of
length~$rn+1$ bordered by $(-\infty,-\infty)$, and
let $R_{rn}^\star$ be the number of $r$--repeated permutations of
length~$rn$ bordered by $(-\infty,+\infty)$
The corresponding  ordinary generating functions 
admit a continued fraction
expansions of the Jacobi type,
\[\begin{array}{l}\ds
\ds \sum_{\nu\ge0} R_{r\nu+1}z^{r\nu+1}
= \\ \qquad \ds
\cfrac{z}{1-2\cdot 1^r\, z^r-\cfrac{1\cdot 2^2\cdots r^2\cdot (r+1)\cdot z^{2r}}
{1-2\cdot (r+1)^r\, z^r-\cfrac{(r+1)\cdot (r+2)^2\cdots (2r)^2\cdot (2r+1)\cdot z^{2r}}
{\ddots}}},\end{array}
\]
\[\begin{array}{l}\ds
\ds \sum_{\nu\ge0} R_{r\nu}^\star z^{r\nu}
= {} \\ \ds \qquad
\cfrac{1}{1-(0^r+ 1^r)z^r-\cfrac{1^2\cdot 2^2\cdots r^2\cdot z^{2r}}
{1-(r^r+(r+1)^r)z^r-\cfrac{(r+1)^2\cdot (r+2)^2\cdots (2r)^2\cdot z^{2r}}
{\ddots}}},
\end{array}
\]
where the numerators are of  degree~$2r$ and the denominators are of
degree~$r$ in the depth index.
\end{proposition}
\begin{proof}
In accordance with Note~\ref{frvi-note} and Figure~\ref{sweepviennot-fig},
it is best to regard the Fran{\c c}on--Viennot
correspondence as the inductive construction of an increasing
binary tree. If at some stage there are $\ell$ dangling links,
then there are several cases to be considered for adding a node depending on
its type. Here is a table giving: the types of elements in a
permutation; the types of nodes in the tree; the number of
possibilities. 
\begin{center}
\begin{tabular}{c|c|c|c|c}
\hline
Perm.: & peak & valley & double rise & double fall \\
\hline
Tree: & leaf & double node & left branching & right branching \\
\hline
Poss.: & $\ell$ & $\ell$ & $\ell$ & $\ell$ \\ 
& $\ell-1$ & $\ell$ & $\ell-1$ & $\ell$ \\
\hline
\end{tabular}
\end{center}
The third line corresponds to allowing the largest (and last inserted)
value $rn+1$ to fall anywhere, the permutation being bordered
by $(-\infty,-\infty)$. The fourth line corresponds to a
permutation to be such that $rn+1$ occurs at  the end, that is, a
permutation bordered by $(-\infty,+\infty)$. The history of the tree
construction corresponds to Motzkin paths that start at altitude 1, 
have steps that are grouped in batches of~$r$, all of the same type
within a batch. This is converted into a standard Motzkin path by a
succession of two operations: $(i)$~shift the path down by~1;
$(ii)$~then divide the altitude by~$r$, so that sequences of steps of a single
type inherit weights multiplicatively. (This is essentially
the argument of~\cite{FlFr89}.) The statement results.
\end{proof}

The close resemblance between the case $r=3$ of
Proposition~\ref{r-prop} and Theorem~\ref{conrad-sec} is striking.
Notice however that it is required to adjust the possibility function
for level steps, i.e., correct the denominators. 
This is achieved by means of ``polarization'':

\begin{definition}
An $r$-repeated permutation of length $3\nu+1$ is said to be polarized  
if \emph{some} (possibly none, possibly all)
of the \emph{consecutive factors} in the word
 representation
of the permutation that are of the form  $3j+3,3j+2,3j+1$ 
or of the form $3j+1,3j+2,3j+3$ are marked.
\end{definition}
We shall use a minus sign as a mark and write $\bar 3 \equiv -3$.
For instance, two polarized 3-repeated permutation of  size~19 are
\[
\hbox{\small\def\ov{\overline}\renewcommand{\arraycolsep}{2truept}
$\begin{array}{cccccccccccccccccccc}
14&2&17&8&19&3&5&6&16&   {12}&   {11}&   {10}&7&13&1&4&18&9&15\\
14&2&17&8&19&3&5&6&16&\ov{12}&\ov{11}&\ov{10}&7&13&1&4&18&9&15
\end{array}$
}
\]
(Only one factor is amenable to marking here.)
We state:
\begin{theorem}
The exponential generating function of 3--repeated polarized permutations 
bordered by $(-\infty,-\infty)$ is
\[
\smh(z).
\]
\end{theorem}
\begin{proof} By Conrad's fractions, Fran{\c c}on--Viennot, and
the reasoning that underlies Proposition~\ref{r-prop}. 
Polarization has been introduced on purpose in order to add to  the number of
possibilities: it contributes
$2(3\ell+1)$ further
possibilities to the Motzkin path when at altitude $\ell$,
corresponding to a (polarized) sequence of three consecutive one-way branching nodes
of two possible types (left- or right-branching) attached to one of the
$3\ell+1$ dangling links.
\end{proof}
The notion of polarization can be similarly introduced to
interpret the coefficients of $\cmh(z)$ (details omitted).
This theorem nicely completes the picture of the relation between
$r$-repeated permutations and special functions. 
\begin{itemize}
\item[$r=1$:] We are dealing with unconstrained permutations. In this
case, Proposition~\ref{r-prop} reduces to the two continued fraction
expansions relative to $\sum n!\, z^n$. Note that marking, with an
additional variable~$t$, the number of rises (i.e., double rises and
valleys), leads to a continued fraction expansion of the bivariate
\emph{ordinary} generating function of Eulerian
numbers~\cite{Flajolet80b}, which was found analytically
Stieltjes~\cite{Stieltjes93}.
\item[$r=2$:] We are dealing with the doubled permutations
of Flajolet-Fran{\c c}on~\cite{FlFr89}. The permutations are
enumerated by Euler numbers, and when rises are taken into account,
one obtains a bivariate ordinary generating function that is the Laplace
transform of the Jacobian elliptic functions $sn(z,t),cn(z,t)$.
\item[$r=3$:] The generating functions are related to $\sm,\cm$, the
correspondence being exact when polarization is introduced. 
\end{itemize}

The Dixonian functions involve once more a third-order
symmetry that is curiously evocative of the fact that they
parametrize the Fermat cubic. It would be of obvious interest (but
probably difficult) to identify which special functions are associated
with higher order symmetries corresponding to $r\ge4$.  

\section{Further connections}\label{connec-sec}

Our goal in this article has been to demonstrate that the Dixonian
parametrization of the Fermat cubic has interesting ramifications 
in several different fields.
The way these functions have largely independently 
surfaced in various domains is striking: occurrences 
now known include
the theory of continued fractions and orthogonal polynomials,
special functions (e.g, Lundberg's hypergoniometric functions),
combinatorial analysis (the elementary combinatorics of permutations),
and a diversity of 
stochastic processes (special urn models and branching processes, 
but also birth and death processes). 
We now briefly discuss other works,
some very recent, which confirm that Dixonian functions should 
indeed be considered as part of the arsenal of special functions.

\subsection{Jacobi and Weierstra{\ss} forms.}\label{jcwei-subsec}

As it is well known, there are three major ways of introducing
elliptic functions~\cite{WhWa27}, namely, by way of the Jacobian 
functions $sn,cn,dn$ (defined from inverses of an Abelian integral over a 
curve $y^2=P_4(x)$), by their Weierstra{\ss} form $\wp$ (associated to a
curve $y^2=P_3(x)$), and by theta functions. 
The reductions of $\sm,\cm$ to normal form are no
surprise since they are granted by general theorems. What stands out,
however, is the simplicity of the connections, which are accompanied in
one case by further combinatorial connections.

\subsubsection{The Jacobian connection.} \label{jacobconnect-sec}
This is the one observed by Cayley in his two page
note~\cite{Cayley82} and already alluded to in
Section~\ref{fermat-sec}. Cayley finds the parametrization of Fermat's
cubic $x^3+y^3=1$ (locally near $(-1,2^{1/3})$)
in the form $(x,y)=(\xi(u),\eta(u))$, where
\[
\xi(u)=\frac{-1+\theta scd}{1+\theta scd},\quad
\eta(u)=\frac{2^{1/3}(1+\theta^2s^2)^2}{1+\theta scd},\quad
\theta:=3^{\frac14}e^{5i\pi/12}
\]
and here $s\equiv sn(u)$, $c\equiv cn(u)$, $d\equiv dn(u)$ are
Jacobian elliptic functions for the modulus $k:=e^{5i\pi/6}$. 
(This corrects what seems to be minor errors in Cayley's
calculations.)
Cayley's calculations imply an expression of $\sm,\cm$ in terms of
$sn,cn,dn$ as simple variants of $\xi,\eta$ (with $\pi_3$ as in~(\ref{pi3def})):
\[
\sm(z)=\xi\left(\frac{z+\pi_3/6}{2^{1/3}\theta}\right),
\qquad
\cm(z)=\eta\left(\frac{z+\pi_3/6}{2^{1/3}\theta}\right).
\]
For an alternative approach based on Lundberg's hypergoniometric functions,
see~\cite[Sec.~5]{LiPe01}.

\subsubsection{The Weierstra{\ss} connection.} The relations appear to
be more transparent than with Jacobian functions.
As a simple illustration, consider the function
\[
P(z):=\smh(z) \cdot \cmh(z).\]
From the basic differential equation framework, this function
satisfies the ordinary differential equations,
\[
P''=6P^2,\qquad P'^2=4P^3+1. 
\]
Thus, up to a shift of the argument, $P(z)$ is a Weierstra{\ss} $\wp$
with parameters $g_2=0$, $g_3=-1$ (corresponding to the usual
hexagonal lattice). While $\wp(0)=\infty$, the initial conditions are
here $P(z)=z+O(z^2)$. In this way we find (see~(\ref{pi3def}))
\[
P(z)=\cal P(z-\zeta_0;0;-1), \qquad
\zeta_0=\frac23\pi_3=\frac{1}{3\pi}\Gamma\left(\frac13\right)^3.
\]
From these calculations, there results that $P=\smh\cdot \cmh$ is
implicitly defined as the solution of the equation
\begin{equation}\label{f131243}
\int_0^Y \frac{dw}{\sqrt{1+4w^3}}\equiv Y \cdot {}_2
F_1\left[\frac13,\frac12;\frac43;-4Y^3\right]=z.
\end{equation}
Thus, the continued fraction expansion of $\sm\cdot \sm$ found in
Section~\ref{conrad-sec} can be re-expressed as follows:
\begin{proposition}
The Laplace transform of the compositional inverse of 
the function $Y \cdot {}_2
F_1\left[\frac13,\frac12;\frac43;-4Y^3\right]$,
equivalently of the special $\wp(z;0,-1)$ expanded near its 
real zero $\zeta_0:=\frac23\pi_3$,
admits a continued fraction
expansion with sextic numerators and cubic denominators.
\end{proposition}
(The hypergeometric parameters differ from those of Proposition~\ref{prop-hyper}.
This function otherwise constitutes a good case of application of Andr\'e's
method of Note~\ref{andre-note}.)

In passing, we remark that $P(z)$ serves to express the
fundamental function (the $\psi$-function) of the 
urn model defined by the matrix 
\[
\cal T_{23}=\left(\begin{array}{cc}-2 & 3\\ 4 & -3\end{array}\right).\]
This model (corresponding to case~$A$ in Figures~\ref{sixell-fig}
and~\ref{sixlat-fig})
is of interest as
it describes the fringe behaviour of 2--3 trees and other locally
balanced trees~\cite{FlGaPe05,PaPr98b}.
The relation is (see~\cite[p.~1223]{FlGaPe05} for the definition
of~$\psi$)
\[
\psi(z)=2^{2/3}P(2^{1/3}z)=2^{2/3}\smh(2^{1/3}z)\cmh(2^{1/3}z),\]
as can be verified either from differential relations or 
from the available $\wp$ forms.
Thus yet another combinatorial model of Dixonian functions
is available in terms of
histories of the $\cal T_{23}$ urn, corresponding to the 
rewrite rules
\[
{\sf xx \longrightarrow yyy}, \qquad
{\sf yyy \longrightarrow xxxx}.\]
In other words:
\begin{proposition}  The ${\cal T_{23}}$ urn model can be described in terms
of Dixonian functions via the product 
$P(z)=\sm(z)\cdot\cm(z)$. In particular, the number of histories of the
$\cal T_{23}$ urn initialized with ${\sf xx}$ such that at time $3\nu+1$ all
balls are of type ${\sf y}$ is 
\[
H_{3\nu+1,3\nu+3}=(3\nu+1)!2^{\nu+1}\, [z^{3\nu+1}] \smh(z)\cdot
\cmh(z).\]
\end{proposition}

In an unpublished manuscript Dumont~\cite{Dumont88}
also observes (by way of singularities)
that the function
\[
Q(z):=\frac{\sm(z)}{3(1-\cm(z))},
\]
which now has a double pole at~0 is none other than $\wp(z;0,\frac{1}{27})$.
From this and similar considerations, Dumont 
then obtains in a simple way the identities
\[
\cm(z)=\frac{3\wp'(z)+1}{3\wp'(z)-1},\qquad
\sm(z)=\frac{6\wp(z)}{1-3\wp'(z)},
\qquad
\wp(z):=\wp(z,0,\frac{1}{27}).
\]

\subsection{Laplace transforms of elliptic functions.}
Transforms of elliptic functions
are only considered sporadically in the literature, usually in
the context of continued fraction theory.
Some of them are explicitly known in the classical case of Jacobian 
elliptic functions, which are somehow close to $\sm,\cm$ given the
developments of Section~\ref{jacobconnect-sec}.

Stieltjes, followed by Rogers, was the first to determine 
continued fraction expansions for the Laplace transforms of the three 
fundamental Jacobian functions, $sn,cn,dn$. Some of these expansions
were rediscovered by Ramanujan who made the  Laplace
transforms explicit. For instance, following Perron's account
(pp.~134--135 and 219 of~\cite{Perron54}) and Berndt's edition of
Ramanujan's Notebooks~\cite[p.~163]{Berndt91}, one has
\begin{equation}\label{stieltram}
\begin{array}{lll}
\ds \int_0^\infty cn(tx;k) e^{-t}\, dt &=& \ds
\frac{2\pi}{kK}\sum_{\nu=0}^\infty
\frac{q^{\nu+\frac12}}{1+q^{2\nu+1}}\,\frac{1}{1+x^2\left(\frac{(2\nu+1)\pi}{2K}\right)^2}
\\
&=&\ds
\cfrac{1}{1+\cfrac{1^2x^2}{1+\cfrac{2^2k^2x^2}{1+\cfrac{3^2x^2}{1+\ddots}}}}.
\end{array}
\end{equation}
(The numerators are of the alternating form
$1^2,2^2k^2,3^2,4^2k^2,5^2x^2,\ldots$; usual notations from elliptic
function theory are employed~\cite{WhWa27}.)
The calculation of the Laplace transform is effected via the Fourier
expansion of $cn$, which goes back to Jacobi and is detailed
in~\cite[\S22.6]{WhWa27}. 
It would be of obvious interest to carry out calculations and
determine in which class
of special functions the Laplace transforms of the Dixonian functions live.

Note otherwise that the Hankel determinant
evaluations stemming from~(\ref{stieltram}) lie at the basis of
Milne's results~\cite{Milne02} regarding sums of squares,
while the continued fraction expansion is related to several
permutation models of Jacobian elliptic functions that were briefly
mentioned in Section~\ref{secondperm-sec}.

\subsection{Orthogonal polynomial systems.} 
As it is well known each continued fraction of type~$J$
is associated to
an orthogonal polynomial system (OPS), which provides
in particular the denominators of the fraction's 
convergents~\cite{Chihara78,Perron54,Wall48}.
The OPS arising from the
Laplace transforms of Jacobian functions have
a long tradition, starting with Stieltjes who
first found an explicit representation of the orthogonality measure.
See the study by Ismail, Valent, and Yoon~\cite{IsVaYo01} for a recent
perspective and pointers to the older literature as well as the book
by Lomont and Brillhart~\cite{LoBr01} entirely dedicated to ``elliptic polynomials''.

It  is a striking fact  that the orthogonal polynomial systems related
to the  continued fraction expansions of Section~\ref{conrad-sec} have
very recently  surfaced   in independent works  of  Gilewicz, Leopold,
Ruffing, and Valent~\cite{GiLeRuVa05,GiLeVa05}.  These authors are 
motivated by the classification  of (continuous-time) birth and  death
processes\footnote{%
	Thanks to  works  of   Karlin and   MacGregor  in the 1950's,   many
	stochastic characteristics of a  process can be described  in
	terms of a continued  fraction, its OPS, and its orthogonality
	measure(s);   see for   instance~\cite{FlGu00}   for a  recent
	review.
}, especially the ones whose
rates are polynomials in the size of the population. 
Here the birth and death rates are cubic.

From the rather remarkable calculations of~\cite{GiLeRuVa05}, one gets
in particular [do $c\to0$ and $\mu_0=0$]
the generating function of a family of polynomials,
\begin{equation}\label{valent}
\begin{array}{lll}
\ds \sum_{n\ge0} Q_n(z)\frac{t^{3n}}{(3n)!} &
=& \ds  (1-t^3)^{-1/3} E_{3,0}\left(z \theta(t)^3\right)\\
E_{3,0}(y) &:=& \ds \sum_{n\ge0} \frac{y^n}{(3n)!},\qquad
\theta(t)\quad:=\quad
\int_0^t \frac{dw}{(1-w^3)^{2/3}}.
\end{array}
\end{equation}
(This OPS resembles a hybrid of the Meixner 
and Brenke classifications of orthogonal polynomials~\cite[Ch.~V]{Chihara78}.)
The first few values, $Q_0,Q_1,Q_2,Q_3,Q_4$, are
\[
\begin{array}{l}
1,\quad 2+z,\quad 160+100\,z+{z}^{2},\quad 62720+42960\,z+672\,{z}^{2}+{z}^{3},
\\
68992000
+49755200\,z+963600\,{z}^{2}+2420\,{z}^{3}+{z}^{4}.
\end{array}
\]
These polynomials are the monic versions of the denominators of the
convergents in the continued fraction
\[
\cfrac{1}{z+2-\cfrac{36}{z+98-\cfrac{14400}{z+572-\ddots}}},
\]
which, up to normalization, represents the Laplace transform of
$\cm(z)$. 
The occurrence in~(\ref{valent})  of the fundamental Abelian integral from which
$\sm,\cm$ are defined is especially striking.
As shown in~\cite{GiLeRuVa05,GiLeVa05}, the moment problem is
indeterminate; see also
Ismail's forthcoming book~\cite{Ismail05} for a perspective.

\subsection*{Conclusion.}
We have studied a class of continued fraction with sextic numerators
and cubic denominators and shown that, thanks to several converging
works, these are endowed with a rich set of properties.
This suggests, regarding the Stieltjes-Ap\'ery fractions described in the
introduction, that it would be of great interest to determine whether 
these are similarly endowed with a rich structure as regards
combinatorics, special functions, and orthogonality relations.
The question of investigating continued fractions whose coefficients
are polynomials of higher degrees is a tantalizing one, but it is
likely to be very difficult.

\bigskip
\noindent
{\bf Acknowledgements}. Eric Conrad
is deeply indebted to Dominique Dumont (Universit\'e Louis Pasteur,
Strasbourg) for 
referring him to the paper of A.\ C.\ Dixon that forms a centerpiece
of  his PhD dissertation~\cite{Conrad02}. Both authors also express their
gratitude to Dumont for many ensuing discussions relative to this
work. We are thankful to Hjalmar Rosengren for first
pointing us to the
works of Lindqvist and Peetre---Peter Lindqvist then very kindly 
offered a wealth of
mathematical and historical material relative to Lundberg's
hypergoniometric functions.
We are finally grateful to Galliano Valent for sharing with us an
unpublished version of his recent work.

\bibliographystyle{amsplain}
\bibliography{algo}

\end{document}